\input amstex
\input epsf
\documentstyle{amsppt}
\magnification=1050
\input xypic
\NoBlackBoxes
\NoRunningHeads
\hsize 6.5truein
\voffset.7in
\document

\define\wt{\widetilde}
\define\ov{\overline}

\define\dfn#1{\definition{\bf\underbar{Definition #1}}}

\define\surj{\twoheadrightarrow}
\define\inj{\rightarrowtail}

\define\prf{\demo{\underbar{Proof}}}
\define\endpf{\enddemo}
\define\tensor{\hat\otimes}
\define\ttensor{\hat\otimes}

\define\exsp{\vskip0.1in}

\define\textit{\it}
\font\script = rsfs10
\define\rel{\Cal R}
\define\RH{\Cal H}

\def\tbrm#1{\text{\bf{#1}}}

\def\aA{\hbox{\script A}}

\def\cC{\hbox{\script C}}

\topmatter
\title
Relatively Hyperbolic Groups, Rapid Decay Algebras and a Generalization of the Bass Conjecture
\endtitle
\vskip.2in
\author
R. Ji (IUPUI), C. Ogle (OSU) and B. Ramsey (IUPUI)
\endauthor
\date
June 2007
\enddate
\email
ronji@math.iupui.edu, ogle@math.mps.ohio-state.edu, bramsey@math.iupui.edu
\endemail
\endtopmatter
\vskip.5in

\proclaim{\bf{Abstract}} By deploying dense subalgebras of $\ell^1(G)$ we generalize the Bass conjecture in terms of  Connes' cyclic homology theory. In particular, we propose a stronger version of the $\ell^1$-Bass Conjecture. We prove that hyperbolic groups relative to finitely many subgroups, each of which posses the polynomial conjugacy-bound property and nilpotent periodicity property, satisfy the $\ell^1$-Stronger-Bass Conjecture. Moreover, we determine the conjugacy-bound for relatively hyperbolic groups and compute the cyclic cohomology of the $\ell^1$-algebra of any discrete group.
\endproclaim

\head Introduction
\endhead
\vskip.2in

Let $R$ be a discrete ring with unit. A classical construction of
Hattori-Stallings ( [Ha], [St] ) yields a homomorphism of abelian groups
$$
K^a_0(R)\overset{Tr^{HS}}\to\longrightarrow HH_0(R)
$$
where $K^a_0(R)$ denotes the zeroth algebraic $K$-group of $R$,
and $HH_0(R) = R/[R,R]$ its zeroth Hochschild homology group.
Precisely, given a projective module $P$ over $R$, one chooses a
free $R$-module $R^n$ containing $P$ as a direct summand,
resulting in a map
$$
End_R(P)\overset{\iota}\to\hookrightarrow End_R(R^n) =
M_{n,n}(R)\overset{Trace}\to\longrightarrow HH_0(R)
$$
One then verifies the equivalence class of $Tr^{HS}([P]) :=
Trace(\iota(Id_P))$ in $HH_0(R)$ depends only on the isomorphism
class of $P$, is invariant under stabilization, and sends direct
sums to sums. Hence it extends uniquely to a map on $K_0(R)$, as
above.

Let $G$ be a discrete group, $k$ a commutative ring and $k[G]$ the
group algebra of $G$ with coefficients in $k\subseteq \Bbb C$. Let
$<G>$ denote the set of conjugacy classes of $G$. Taking tensor
products over $k$, $HH_0(k[G]) = \oplus_{x\in <G>}HH_0(k[G])_x =
\oplus_{x\in <G>}k$ decomposes as a direct sum of copies of $k$
indexed on $<G>$. Let $\pi_x : HH_0(k[G])\to k$ be the projection
which is the identity on the summand indexed by $x$ and sends the
other summands to zero. Given a finitely-generated projective
module $P$ over $k[G]$, the $P$-rank of an element $g\in G$ is
given by $r_P(g) := \pi_{<g>}Tr^{HS}([P])$. By convention, the
conjugacy classes associated to elements of finite order are
referred to as \underbar{elliptic}, and the sum of all such the
\underbar{elliptic summand}. In [Ba], H. Bass stated the
following conjecture.

\proclaim{\bf\underbar{Bass Conjecture}} If $P$ is a
finitely-generated projective module over $\Bbb Z[G]$, then
$r_P(g) = 0$ if $g\ne id$.
\endproclaim

Although still unknown in general, many partial results are known,
beginning with the case $G$ a linear group, verified by Bass
himself in [Ba]. For the integral group algebra, Linell has shown that $r_P(g) = 0$
whenever $1\ne g$ is of finite order. Also, the inclusion
$\Bbb Z\hookrightarrow \Bbb C$ induces an inclusion $HH_0(\Bbb
Z[G]) = \oplus_{x\in <G>}\Bbb Z\hookrightarrow HH_0(\Bbb C[G]) =
\oplus_{x\in <G>}\Bbb C$ (where $HH_0(\Bbb C[G])$ is computed with
tensors over $\Bbb C$) which is the obvious inclusion on each
summand. Thus the Bass conjecture is implied by the

\proclaim{\bf\underbar{Strong Bass Conjecture - SBC}} For each non-elliptic class $x$, the image of
the composition $\pi_x\circ Tr^{HS}:K^a_0(\Bbb C[G])\to HH_0(\Bbb C[G])\surj HH_0(\Bbb C[G])_x$ is zero.
\endproclaim

Henceforth all algebras and tensor products are assumed to be over
$\Bbb C$. If $A$ is a Fr\'echet algebra and $A^{\delta}$ denotes
$A$ with the discrete topology, then there is a commuting diagram
$$
\diagram
K_0^a(A^{\delta})\dto\rto^{Tr^{HS}} & HH_0(A^{\delta})\dto\\
K_0^t(A)\rto^{ch_0} & HC^t_0(A)
\enddiagram
$$
where the vertical maps are induced by the continuous map
$A^{\delta}\to A$ which is the identity on elements, $ch_0$ is the
Connes-Karoubi-Chern character in dimension zero [CK], and
$HC^t_0(A)$ the zeroth (reduced topological) cyclic homology group
of $A$. By [Ti], this same diagram exists for arbitrary topological algebras $A$
when $A$ is topologized
by the fine topology, in which case one has an equality $HC^a_*(A)
= HC^t_*(A)$. The advantage to working with $ch_*$ is
that it fits into a family of Chern characters
$$
ch_n^m: K_n^t(A)\to HC^t_{n + 2m}(A)
$$
which satisfy the identity $S\circ ch_n^m = ch_n^{m-1}$ for $m\ge
1$, where $S: HC^t_{n+2m}(A)\to HC^t_{n+2(m-1)}(A)$ is the $S$
operator. Moreover, For when $A$ is a (locally convex) algebra topologized 
with the fine topology, there is an isomorphism $ HC_*(A) = HC^t_*(A)$, with the left-hand side
denoting \underbar{algebraic} cyclic homology. In particular, we can take $A = \Bbb C[G]$ equipped
with the fine topology. Then for each conjugacy class $x$, there is a summand 
$HC_*(\Bbb C[G])_{x}$ of $HC_*(\Bbb C[G]) = HC_*^t(\Bbb C[G])$, and a corresponding projection
operator $(\pi_x)_*: HC_*(\Bbb C[G])\surj HC_*(\Bbb C[G])_x$.
In light of this, it is natural to consider the following stronger version of
SBC.

\proclaim{\bf\underbar{Stronger Bass conjecture - SrBC}} For each non-elliptic class $x$,
the image
of the composition $(\pi_x)_*\circ ch_*:K^t_*(\Bbb C[G])\to HC_*(\Bbb C[G])\to HC_*(\Bbb C[G])_x$ is zero.
\endproclaim

It follows by [Og1] that the image of $ch_*:K^t_*(\Bbb
C[G])\to HC_*(\Bbb C[G])$ is \underbar{at least} the elliptic
summand; the non-trivial part of the conjecture is
to show it is no more than this.

For a given non-elliptic class
$x$, we say that \underbar{$x$ satisfies the
nilpotency condition} if the operator $S^*:HC^*(\Bbb C[G])_x\to HC^{*+2}(\Bbb
C[G])_x$ is nilpotent. We say that \underbar{$G$ satisfies the
nilpotency condition} if $x$ does for each non-elliptic conjugacy
class $x$. Combining the results of Connes-Karoubi with those of
Burghelea, one has the the following immediate consequence, used
in [Ec1], [Ji3] and [Em]:

\proclaim{\bf\underbar{Observation}} Let $x$ be a non-elliptic
conjugacy class satisfying the nilpotency condition.  Then the
composition $K_n(\Bbb C[G])\to HC_n(\Bbb
C[G])\overset{\pi_x}\to\longrightarrow HC_n(\Bbb C[G])_x$ is zero
for all $n\ge 0$. Consequently, if $G$ satisfies the nilpotency
condition, then the SrBC is true for $G$.
\endproclaim

Unlike SrBC, there are known examples of discrete groups which do
not satisfy the nilpotency condition. However, there is a large
class of groups which do. The most inclusive results are those of
Emmanouil  [Em]. Following (loc. sit.), we denote by $\cC(\Bbb
C)$ the class of groups which satisfy the nilpotency condition
over $\Bbb C$ (this is the same as over $\Bbb Q$, but the notation
is more in keeping with the rest of our paper). \vskip.1in

Let $\ell^1(G)$ denote the $\ell^1$ algebra of $G$. Fixing a
proper word-length function $L$ on $G$, we denote by
$H^{1,\infty}_L(G)$ the standard $\ell^1$ rapid-decay algebra of
$G$; this is a Fr\'echet subalgebra of $\ell^1(G)$ containing
$\Bbb C[G]$ which is smooth in $\ell^1(G)$ - that is, the
inclusion $H^{1,\infty}_L(G)\hookrightarrow \ell^1(G)$ induces an
isomorphism on topological $K$-groups in all degrees [Jo1], [Jo2]. The
Connes-Karoubi-Chern character therefore produces a map
$$
K^t_*(\ell^1(G))\overset\cong\to\longleftarrow
K^t_*(H^{1,\infty}_L(G))\to HC^t_*(H^{1,\infty}_L(G))
$$
Analogous to the case of the group algebra, there is a naturally defined collection
of summands $\{HH^t_*(H^{1,\infty}_L(G))_x\}$ resp. $\{HC^t_*(H^{1,\infty}_L(G))_x\}$
of $HH^t_*(H^{1,\infty}_L(G))$ resp. $HC^t_*(H^{1,\infty}_L(G))$ indexed on the set of conjugacy classes of $G$, with a similar result for topological Hochschild and cyclic  cohomology. This leads to
the following variant of SrBC, which was one of our original motivations for this paper.

\proclaim{\bf\underbar{$\ell^1$-SrBC}} For each non-elliptic conjugacy class $x$, the image of
the composition $(\pi_x)_*\circ ch_*:K^t_*(\ell^1(G))\to HC^t_*(H^{1,\infty}_L(G))\surj 
HC^t_*(H^{1,\infty}_L(G))_x$ is zero.
\endproclaim

\bigskip

The $\ell^1$-Bass Conjecture states that the $\ell^1$
Hattori-Stallings trace: $HS^1: K_0(\ell^1(G))\longrightarrow HH_0(\ell^1(G))$ 
has a range concentrated in elliptic summands. In [BCM] it is shown that the Baum-Connes 
Conjecture for a discrete group implies its $\ell^1$-Bass Conjecture. Thus, a-T-menable groups, 
in particular, amenable groups satisfy the $\ell^1$-Bass conjecture, following the work of 
Higson and Kasporov [HK]. The $\ell^1$-SrBC we have proposed here is stronger than
the $\ell^1$-Bass Conjecture. Our method requires determining the nilpotency of the periodicity operator on cyclic cohomology and some geometric analysis on the group. However, it does not require the verification of Baum-Connes conjecture.
\bigskip

An outline of the paper is as follows. In section 1, following a
preliminary review of Hochschild, cyclic and periodic cyclic
(co)homology, we introduce in section 1.2 the notion of a
\underbar{bounding function}. A family $\Cal B$  of such
functions, together with a word-length function $L$ on $G$, equips
$\Bbb C[G]$ with a collection of semi-norms, for which the
completion $H_{\Cal B,L}(G)$ of $\Bbb C[G]$ is a subalgebra of
$\ell^1(G)$. In this paper the term \underbar{rapid decay algebra}
refers to \underbar{any} subalgebra of $\ell^1(G)$ of this type.
In fact, the rapid decay algebra $H^{1,\infty}_L(G)$ occurs as
$H_{\Cal P,L}(G)$, where $\Cal P$ is the family of polynomial
bounding functions. The category $\underline{\Cal B}$ of all
families of bounding functions is a poset, with partial ordering
given by inclusion. This  category has an initial object $\Cal
B_{min}$, a terminal object $\Cal B_{max}$, and for all $(G,L)$
(meaning a discrete group $G$ equipped with a word-length function
$L$) there are equalities $H_{\Cal B_{min},L}(G) = \ell^1(G)$,
$H_{\Cal B_{max},L}(G) = \Bbb C[G]$. Moreover, for each $(G,L)$,
the association $\Cal B\mapsto H_{\Cal B,L}(G)$ defines a
contravariant functor $F_{G,L}:\underline{\Cal B}\to (top.\ alg.)$
such that the unique inclusion $\Cal B_{min}\inj \Cal B_{max}$
maps under $F_{G,L}$ to the natural inclusion $\Bbb C[G]\inj
\ell^1(G)$.  The result is a class of rapid decay algebras
indexed by $\underline{\Cal B}$, lying between the group algebra $\Bbb C[G]$
and the $\ell^1$ algebra $\ell^1(G)$. Now there is a
\underbar{topological} as well as  a \underbar{bornological}
approach to extending (co)homology theories (Hochschild, cylcic,
periodic cyclic, bar,etc\dots) to non-discrete algebras. In
general, these extensions apply to different categories and hence
yield different theories. However, for Fr\'echet algebras these
extensions are the same [Me2]. For $\Cal B$
countable, $H_{\Cal B,L}(G)$ is a Fr\'echet subalgebra of
$\ell^1(G)$; consequently, on the subclass of rapid decay algebras
$\{H_{\Cal B,L}(G)\ |\ \Cal B \text{ countable}\}$, these two
extensions agree. This is formalized in section 1.3. Consequently,
we write $F_*^t(H_{\Cal B,L}(G))$ resp. $F^*_t(H_{\Cal
B,L}(G))$ to denote either topological or bornological homology
resp. cohomology groups, where $F$ represents Hochschild, cyclic,
periodic cyclic or bar homology. Moreover, as we further show in
this section, the cohomology theories identify with the the
corresponding \underbar{$\Cal B$-bounded} theory computed as the
cohomology of the subcocomplex of $\Cal B$-bounded cochains in the
cocomplex associated to the group algebra. This idea, which is a
natural extension of both [Og2] and [M1], is most conveniently cast
in the context of $\Cal B$-bounded cohomology theories associated
to weighted simplicial sets (also defined in section 1.3). To
illustrate, there is an isomorphism
$$
HH^*_t(H_{\Cal B,L}(G))\cong HH^*_{\Cal B}(\Bbb C[G])
$$
where the groups on the left are the topological Hochschild
cohomology groups of $H_{\Cal B,L}(G)$, while the groups on the
right are the $\Cal B$-bounded Hochschild cohomology groups of
$\Bbb C[G]$, computed as the $\Cal B$-bounded cohomology of the
weighted simplicial set $N.^{cy}(G) =$ the cyclic bar construction
on $G$, with weight function induced by $L$. Additionally, these
groups contain summands indexed by the conjugacy classes of $G$, and which
may be identified with appropriate $\Cal B$-bounded
cohomology groups of weighted simplicial summands of $N.^{cy}(G)$.
This is made precise in section 1.4. To further identify the
summands in the manner suggested by Burghelea's result, we need to
consider various (albeit natural) restrictions on $G$. Under the
assumption that $G$ has  a \underbar{$\Cal B$-solvable conjugacy-bound}
there is for each conjugacy class $x\in <G>$ an
isomorphism $HH^*_{\Cal B}(\Bbb C[G])_x\cong H^*_{\Cal
B}(BG_x;\Bbb C)$, where $G_x$ denotes the centralizer of $h$ in
$G$ and $<h> = x$ (Cor. 1.4.6). In the language of group
theory, $\Cal B$-solvability of the conjugacy-bound is
essentially saying that for each conjugacy class there is a
bounding function $\phi\in\Cal B$ for which the length of a conjugator $g$
for two conjugate elements $u$ and $v$ is bounded by $\phi$ of the lengths 
of $u$ and $v$. Again, for each conjugacy class $x$
$F_*^t(H_{\Cal B,L}(G))$ contains a naturally defined summand $F_*^t(H_{\Cal B,L}(G))_x$
($F = HH,HC, HPer$), which suggests the following extension of $\ell^1$-SrBC.

\proclaim{\bf\underbar{$\Cal B$-SrBC}} For each non-elliptic class $x$, the image of
the composition 
$(\pi_x)_*\circ ch_*:K^t_*(H_{\Cal B,L}(G))\to HC_*^t(H_{\Cal B,L}(G))\surj HC_*^t(H_{\Cal B,L}(G))_x$ is
zero.
\endproclaim

We say $G$ satisfies the \underbar{$\Cal B$-nilpotency condition}
if the operator $S^*_x:HC^*(H_{\Cal B,L}(G))_x\to HC^{*+2}(H_{\Cal
B,L}(G))_x$ is nilpotent for each non-elliptic conjugacy class
$x$. As above, it follows by Connes-Karoubi that

\proclaim{\bf\underbar{Observation}} If $G$ satisfies the $\Cal
B$-nilpotency condition, then the conjecture $\Cal B$-SrBC is true
for $G$.
\endproclaim

In section 2, we use the machinery of section 1 to verify the $\Cal B$-nilpotency
condition in a number of interesting cases.  Although not a theorem, in practice it
is generally the case that if $G$ satisfies the $\Cal
B$-nilpotency condition and $\Cal B\subset \Cal B'$, then $G$
satisfies the $\Cal B'$-nilpotency condition. The class of groups
satisfying the $\Cal B$-nilpotency condition will be denoted by
$\cC_{\Cal B}(\Bbb C)$. The weakest condition is $\Cal
B_{max}$-nilpotency, which is nothing else but the nilpotency
condition defined above for the group algebra $\Bbb C[G]$ (i.e.,
$\cC_{\Cal B_{max}}(\Bbb C) = \cC(\Bbb C)$). As a warm-up, and to
illustrate this method, we show in section 2.1 that $G$ satisifes
$\Cal B$-nilpotency for all $\Cal B\ne \Cal B_{min}$ when $G$ is
(i) f.g. nilpotent or (ii) word-hyperbolic. This includes the case
$\Cal B = \Cal P$, which are results due to the first author
[Ji1], [Ji2]. In section 2.2 we extend this by showing

\proclaim{\bf\underbar{Theorem A}} If $\Cal P\subset \Cal B$ and
$G$ is hyperbolic relative to a finite set of subgroups $H_i$
where each $H_i$ is in the class $\cC_{\Cal B}(\Bbb C)$ and has a 
$\Cal B$-solvable conjugacy bound, then $G$
is in the class $\cC_{\Cal B}(\Bbb C)$. In particular, If $\Cal B = \Cal P$, 
then $G$ satisfies the $\ell^1-$SrBC.
\endproclaim

Theorem A depends heavily on the following theorem on relatively hyperbolic groups:

\proclaim{\bf\underbar{Theorem B}} Let $G$ be a finitely
generated group with generating set $X$, such that $G$ is
relatively hyperbolic with respect to $\{ H_\lambda \}_{\lambda
\in \Lambda}$, each member of which has $\Cal P$-solvable conjugacy-bound property.
Then $G$ has a $\Cal P$-solvable conjugacy-bound.  
\endproclaim

The conjugacy problem for relatively hyperbolic groups was solved in [Bu]. However, the $\Cal P$-solvable conjugacy-bound property described in Theorem B is new. Recall that the \underbar{idempotent conjecture} for $\ell^1(G)$ states that the only idempotent elements
 in $\ell^1(G)$ are $0$ and $1$  if $G$ is torsion-free. By a standard argument we have

\proclaim{\bf\underbar{Corollary C}} If $G$ satisfies the hypothesis of Theorem B and is torsion-free, then the idempotent conjecture holds for $\ell^1(G)$.
\endproclaim

\bigskip

In section 2.3, we consider the case when $G_x$ is synchonously
combable for non-elliptic $x$; the results of this section are
preliminary, and indicate a need for further study. One important family of
combable groups are those which are
semi-hyperbolic in the sense of [AB]. As we have already seen,
$\Cal B$-solvability is needed at the non-elliptic classes in
order to properly identify the resulting summands of the
topological Hochschild and cyclic groups. Our second main result
is

\proclaim{\bf\underbar{Theorem D}} If $G$ is semi-hyperbolic, and
$\Cal B$ is any class containing $\Cal E$ (the exponential
bounding class), then $G$ is in $\cC_{\Cal B}(\Bbb C)$ if it is in
$\cC_{\Cal B_{max}}(\Bbb C)$. \endproclaim

For certain classes of semi-hyperbolic groups, this result can be
strengthened to where $\Cal E$ is replaced by $\Cal P$ in the
above theorem. Finally in section 2.4, we deal with the class $\Cal
B_{min}$ (corresponding to the full $\ell^1$-algebra $\ell^1(G))$.
In marked contrast to the structure of the topological cyclic homology
groups of $H_{\Cal B,L}(G)$ when $\Cal B\ne \Cal B_{min}$, we have

\proclaim{\bf\underbar{Theorem E}} For each conjugacy class $x$,
there is an isomorphism
$$
HC^*_t(\ell^1(G))_x\cong H^*_b(BG_x;\Bbb C)\otimes HC^*(\Bbb C)
$$
where the $S^*_x$ map on the left identifies with $Id\otimes S^*$
on the right, and $H^*_b(_-)$ denotes bounded cohomology. Thus in
the case of the $\ell^1$-algebra, the map $S^*_x$ is never
nilpotent for any conjugacy class $x$. Consequently, the class
$\cC_{\Cal B_{min}}(\Bbb C)$ is empty.
\endproclaim

The essential reason for this result is the vanishing of bounded
cohomology for amenable groups in positive degrees. It is this
theorem that motivates consideration of larger classes of bounding
functions (or equivalently, the restriction to rapid decay
subalgebras of $\ell^1(G)$). \vskip.1in

Finally, it is natural to expect a connection between the conjecture $\Cal B$-SrBC and the 
corresponding Baum-Connes conjecture for the rapid decay algebra $H_{\Cal B,L}(G)$.
In fact, the rational surjectivity of the (appropriately defined) 
Baum-Connes assembly map for $K_*^t(H_{\Cal B,L}(G))$ implies $\Cal B$-SrBC
(a very special case of this result was shown in [BCM]). A proof of this
result is provided in the appendix.\vskip.4in

Some of the ideas of the paper grow out of the workshop \lq\lq The property of rapid deca\rq\rq in Jan. 2006 at the American Institute of Mathematics. The first and third authors wish to thank the institute for providing active and friendly atmosphere for the workshop. 
\vskip.3in

\head 1. Hochschild and cyclic (co)homology of general rapid decay
algebras
\endhead
\vskip.2in

\subhead 1.1. Preliminaries
\endsubhead
\vskip.1in

We recall the construction of Hochschild and cyclic homology
([Co1], [Lo]), assuming throughout that the base field is $\Bbb C$,
with all algebras and tensor products being over $\Bbb C$. For an
associative unital algebra $A$ and an $A$-bimodule $M$, let
$C_n(A,M) = M\otimes A^{\otimes n}$. Define $b_n : C_n(A,M)\to
C_{n-1}(A,M)$ by the equation
$$
b_n(m,a_1,\dots, a_n) = (ma_1,a_2,\dots,a_n) +
\sum_{i=1}^{n-1} (-1)^i (m,\dots, a_ia_{i+1},\dots, a_n)
+(-1)^n (a_nm, a_1,\dots, a_{n-1})
$$
This is a differential, and the resulting complex
$C_*(A,M) = (C_*(A,M),b)$ is the Hochschild complex of $A$ with coefficients on $M$.
Its homology is denoted $HH_*(A,M)$. When $M = A$, viewed in the obvious way as
an $A$-bimodule, we write $C_*(A,A)$ simply as $C_*(A)$, referred to as the Hochschild
complex of $A$. The complex $C_*(A)$ admits a degree one chain map $B_*$ given by
$$
\split
B_n(a_0,a_1,\dots,a_n) =
&\sum_{i=o}^n (-1)^{ni}(1, a_i, a_{i+1},\dots, a_n, a_0, a_1\dots, a_{i-1})\\
& -\sum_{i=0}^n(-1)^{ni}(a_i,1,a_{i+1},\dots,a_n,a_0,\dots,a_{i-1})
\endsplit
$$
Let $\Cal B'(A)_{**}$ denote the first-quadrant bicomplex  with
$\Cal B'(A)_{2p,q} = \otimes^{q+1} A$ for $p\ge 0$, $\Cal
B'(A)_{2p,q} = 0$ for $p < 0$ and $\Cal B'(A)_{2p-1,q} = 0$ for
all $p$. The differentials are $b_q:\Cal B'(A)_{2p,q}\to \Cal
B'(A)_{2p,q-1}$ and $B_q:\Cal B'(A)_{2p,q}\to \Cal
B'(A)_{2p-2,q+1}$ for $p\ge 1$. The associated total complex is
the cyclic complex of $A$, denoted $CC_*(A)$; because we are over
a field of characteristic $0$, this is quasi-isomorphic to the
smaller complex $(C_*(A)/(1-\tau_*),b_*)$ (which we also denote by
$CC_*(A)$ when there is no confusion), where
$\tau_n(a_0,\dots,a_n) = (-1)^n(a_n,a_0,\dots,a_{n-1})$. One may
also form the larger 2-periodic bicomplex $\Cal B(A)_{**}$ with
$\Cal B(A)_{2p,q} = \otimes^{q+1} A$, $\Cal B(A)_{2p-1,q} = 0$.
The differentials are $b_q:\Cal B(A)_{2p,q}\to \Cal B(A)_{2p,q-1}$
and $B_q:\Cal B(A)_{2p,q}\to \Cal B(A)_{2p-2,q+1}$, as before. The
total complex $T_*(A) = (\underset{2p+q = *}\to\oplus \Cal
B(A)_{2p,q}, b + B)$ of the resulting bicomplex $\Cal B(A)_{**} =
\{\Cal B(A)_{**}; b,B\}$ is clearly $\Bbb Z/2$-graded (and
contractible). There is a decreasing filtration by subcomplexes
$F^kT_*(A) = (b(C_k(A))\underset n\to\bigoplus C_{n\ge k}(A), b +
B)$, and the periodic cyclic complex $CC_*^{per}(A)$ is the
completion of $T_*(A)$ with respect to this filtraton:
$CC_*^{per}(A) = \underset n\to\varprojlim T_*(A)/F^kT_*(A)$, with
homology $HC_*^{per}(A)$. The dual constructions yield Hochschild,
cyclic and periodic cyclic cohomology. Finally we recall the definition of
Hochschild homology and cohomology in terms of derived functors. For a
unital algebra $A$, $A^{op}$ denotes $A$ with the opposite multiplication.
An $A$-bimodule $M$ is naturally an $A\otimes A^{op}$-module via the
operation $(A\otimes A^{op})\otimes M\ni (a\otimes b)\otimes m\mapsto amb$.
In particular, $A$, via its canonical $A$-bimodule structure, is an
$A\otimes A^{op}$-module. One then has $HH_*(A,M) = Tor_*^{A\otimes A^{op}}(A,M)$,
$HH^*(A,M) = Ext^*_{A\otimes A^{op}}(A,M)$.
\vskip.2in

{\bf\underbar{Note 1.1.1}} The Hochschild \underbar{cohomology} of
$A$, denoted $HH^*(A)$, is the cohomology of the cocomplex
$(Hom(C_*(A),\Bbb C), b^*)$. In terms of the above, this
identifies with $HH^*(A,A^*)$, where $A^*$ is the $A$-bimodule of
linear functionals on $A$. Similarly, the cyclic resp. periodic
cyclic cohomology of $A$ is computed as $HC^*(A) = H^*(CC^*(A))$
resp. $HPer^*(A) = H^*(CPer^*(A))$, where $CC^*(A) :=
Hom(CC_*(A),\Bbb C)$ resp. $CPer^*(A) := Hom(CPer_*(A),\Bbb C)$.
\vskip.2in

A fundamental result and starting point for the computations in
this paper is the well-known computation of the cyclic homology of
$\Bbb C[G]$, due to Burghelea [Bur], which we summarize as:

\proclaim{\bf\underbar{Theorem 1.1.2 -- Burghelea}} \vskip.1in

(i) There are canonical isomorphisms
$$
\gather
HH_*(\Bbb C[G])\cong \oplus_{x\in <G>}HH_*(\Bbb C[G])_x\\
HC_*(\Bbb C[G])\cong \oplus_{x\in <G>}HC_*(\Bbb C[G])_x
\endgather
$$
with the maps in the Connes-Gysin long-exact preserving this decomposition.
\vskip.05in

(ii) For all $x$, there is a canonical isomorphism $HH_*(\Bbb
C[G])_x = H_*(BC_g;\Bbb C)$, where $C_g$ denotes the centralizer
of $g\in G$ ($<g> = x$). \vskip.05in

(iii) If $x$ is elliptic, then $HC_*(\Bbb C[G])_x = H_*(BC_g;\Bbb
C)\otimes HC_*(\Bbb C)$, where $S_x: HC_*(\Bbb C[G])_x\to
HC_{*-2}(\Bbb C[G])_x$ is given by $S_x = Id\otimes S_*^{\Bbb C}$,
$S_*^{\Bbb C}: HC_*(\Bbb C)\to HC_{*-2}(\Bbb C)$. \vskip.05in

(iv) If $x$ is non-elliptic, then $HC_*(\Bbb C[G])_x =
H_*(BN_g;\Bbb C)$, where $N_g = C_g/(g) =$ the quotient of $C_g$
by the infinite cyclic subgroup generated by $g$. Under this
identification, the summand of the Connes-Gysin sequence indexed
by $x$ identifies with the Gysin sequence in homology (with
coeff.s in $\Bbb C$) associated to the fibration $S^1 = \Bbb Z\to
BC_g\to BN_g$. \vskip.05in

(v) The analogues of (i) - (iv) hold for the Hochschild and cyclic
cohomology of $\Bbb C[G]$ (with direct sum replaced by direct
product in (i)).
\endproclaim
\vskip.2in

\subhead 1.2. Bounding functions and rapid decay algebras
\endsubhead
\vskip.1in

We begin with

\proclaim{\bf\underbar{Definition 1.2.1}} A \underbar{bounding
class} $\Cal B$ is a non-empty set of non-decreasing functions
$\{f : \Bbb R_+\to \Bbb R^+\}$ satisfying the following conditions
\vskip.05in

(i) It is closed under the operation of taking positive rational linear combinations.
\vskip.05in

(ii) If $f_1,f_2\in\Cal B$, there is an $f_3\in\Cal B$ such that $f_1\circ f_2\le f_3$.
\vskip.05in

(iii) it contains the constant function $1$.
\endproclaim

By (1.2.1)(iii) The smallest bounding class is $\Cal B_{min} =
\Bbb Q^+ = $ the set of constant functions with value a positive
rational number. There is also a largest class $\Cal B_{max} =
\{f:\Bbb R_+\to \Bbb R^+\ |\ f\text{ is non-decreasing}\}$. Other
basic classes that will appear in this paper are $\Cal L := \{a +
bx\ |\ a,b\in\Bbb Q^+\}$, $\Cal P := $ the set of positive
rational linear combinations of $\{(1+x)^m\}_{m\in \Bbb N}$, and
$\Cal E :=$ the closure, under composition, of the set of positive rational linear combinations of
$\{C^x\}_{C\in\Bbb Q,C \ge 1}$. Two bounding classes $\Cal B$ and $\Cal B'$ are \underbar{equivalent}, written $\Cal B\simeq \Cal B'$, if for each $f$ in one class, there is an $f'$ in the other with $f\le f'$. For two bounding classes $\Cal B$, $\Cal B'$ we write $\Cal B + \Cal B'$ to indicate the bounding class generated by $\Cal B\cup \Cal B'$; this is the smallest bounding class
containing both $\Cal B$ and $\Cal B'$. If $\Cal B\simeq\Cal B +\Cal B'$, we write
$\Cal B'\preceq \Cal B$, and $\Cal B' \prec \Cal B$ if $\Cal B'\preceq \Cal B$ but $\Cal B\not\preceq \Cal B'$. Thus $\Cal B_{min}\prec \Cal L\prec \Cal P\prec \Cal E\prec \Cal B_{max}$.\vskip.2in

A \underbar{weighted} set $(X,w)$ will refer to a countable
discrete set $X$ together with a function $w:X\to\Bbb R_+$. Let
$C_c(X)$ denote the space of complex-valued functions on $X$ with
compact (i.e., finite) support. A morphism $\phi:(X,w)\to (X',w')$
is \underbar{$\Cal B$-bounded} if
$$
\forall f\in\Cal B \ \exists f'\in\Cal B\text{ such that }f(w'(\phi(x)))\le
f'(w(x))\qquad\forall x\in X
$$
Given a weighted set $(X,w)$ and a bounding class $\Cal B$, one
may form the topological space
$$
H_{\Cal B,w}(X) := \left\{\phi:X\to \Bbb C\ \left|\ |\phi|_f :=
\sum_{x\in X} |\phi(x)|f(w(x)) < \infty\right.\right\}
\tag1.2.2
$$
A \underbar{$\Cal B$-isomorphism} $\phi:(X,w)\to (X',w')$ is an
isomorphism of sets which is $\Cal B$-bounded, and which admits a
two-sided inverse which is also $\Cal B$-bounded. There is a natural embedding $H_{\Cal B,w}(X)\hookrightarrow \ell^1(X)$, and so we will typically represent elements of $H_{\Cal B,w}(X)$ in terms of their corresponding (infinite) sums $H_{\Cal B,w}(X)\ni\phi\Leftrightarrow \sum_{x\in X} \phi(x) x$.

In most applications below, $\Cal B$ is countable, in which case
$H_{\Cal B,w}(X)$ is a Fr\'echet space which alternatively may be
viewed as the completion of $C_c(X)$ with respect to the
collection of semi-norms $\{|_-|_f\}_{f\in\Cal B}$.

 \proclaim{\bf\underbar{Proposition 1.2.3}} Let $(X,w)$, $(X',w')$
 be weighted sets, $\Cal B$, $\Cal B'$ bounding classes and $(G,L)$
 a discrete group with proper word-length function.
 \vskip.05in

(i) $H_{\Cal B_{min},w}(X) = \ell^1(X)$. In particular, $H_{\Cal
B_{min},L}(G) = \ell^1(G)$, the $\ell^1$-algebra of $G$.
\vskip.05in

(ii) $H_{\Cal B_{max},w}(X) = C_c(X)$. In particular, $H_{\Cal
B_{max},L}(G) = \Bbb C[G]$, the complex group algebra of $G$.
\vskip.05in

(iii) If $\Cal B\simeq \Cal B'$ then $H_{\Cal B,w}(X) = H_{\Cal B',w}(X)$.
\vskip.05in

(iv) If $\Cal B'\prec \Cal B$, there is a natural inclusion
$H_{\Cal B,w}(X)\subseteq H_{\Cal B',w}(X)$. \vskip.05in

(v) If $\phi:(X,w)\to (X',w')$ is $\Cal B$-bounded, it induces a
continuous map of topological vector spaces $\phi : H_{\Cal
B,w}(X)\to H_{\Cal B,w'}(X')$. If $\phi$ is a $\Cal
B$-isomorphism, this induced map is a continuous isomorphism of
topological vector spaces. \vskip.05in

(vi) $H_{\Cal B,w}(X)\cap H_{\Cal B',w}(X) = H_{\Cal B+\Cal B',w}(X)$.
\vskip.05in

(vii) For all $(G,L)$, there is an equality $H_{\Cal P,w}(G) =
H^{1,\infty}_L(G)=$ the $\ell^1$-Schwartz algebra of $G$ (this
algebra is the completion of $\Bbb C[G]$ in the seminorms
$\{\nu_k\}_{k\in\Bbb N}$ where $\nu_k(\lambda) := \sum_{g\in G}
|\lambda(g)|(1+L(g))^k$; cf. [Jo1], [Jo2]). \vskip.05in

(viii) For all $(G,L)$ and bounding class $\Cal B$, $H_{\Cal
B,L}(G)$ is a subalgebra of $\ell^1(G)$.
\endproclaim

\prf The only point that is not immediate is $(viii)$. To verify
this last property, it suffices to show $H_{\Cal B,L}(G)$ is
closed under multiplication. Given $f\in\Cal B$, we have
$$
\gather
\left|\left(\sum\lambda_{g_1} g_1\right)\left(\sum\lambda_{g_2} g_2\right)\right|_f\\
= \sum_g\left|\sum_{g_1g_2 = g}\lambda_{g_1}\lambda_{g_2}\right|f(L(g))\\
\le \sum_g\left(\sum_{g_1g_2 = g}\left|\lambda_{g_1}\lambda_{g_2}\right|f(L(g_1)+L(g_2))\right)\\
\le \sum_g\left(\sum_{g_1g_2 = g,L(g_1)\le L(g_2)}\left|\lambda_{g_1}\lambda_{g_2}\right|f(2L(g_2))\right) +
\sum_g\left(\sum_{g_1g_2 = g,L(g_2)\le L(g_1)}\left|\lambda_{g_1}\lambda_{g_2}\right|f(2L(g_1))\right)\\
\le \left|\sum\lambda_{g_1} g_1\right|_1\left|\sum\lambda_{g_2} g_2\right|_{f_2} +
\left|\sum\lambda_{g_1} g_1\right|_{f_2}\left|\sum\lambda_{g_2} g_2\right|_1< \infty
\endgather
$$
where $\left|\sum\lambda_{g} g\right|_1$ is the $\ell^1$-norm and
$f_2$ denotes any element of $\Cal B$ satisfying $f(2x)\le f_2(x)$
for all $x$ (this element exists by (1.2.1) (ii)).\hfill //
\endpf

{\bf\underbar{Note}} The term \underbar{rapid decay subalgebra of
$\ell^1(G)$} conventionally refers to the Schwartz algebra
$H^{1,\infty}_L(G)$. For this paper, a rapid decay subalgebra of
$\ell^1(G)$ will mean any one of the subalgebras $H_{\Cal B,L}(G)$
associated to a bounding class $\Cal B$ and a group with
word-length $(G,L)$.

\proclaim{\bf\underbar{Lemma 1.2.4}} Let $(X,w)$ be as above, with
$w$ a proper weight function, and $\Cal B$ a countable bounding
set. Let $H_{\Cal B,w}(X)^*$ denote the space of continuous linear
functionals on $H_{\Cal B,w}(X)$, and $Hom_{\Cal B}(X,\Bbb C)$ the
vector space of $\Cal B$-bounded functions $f:X\to\Bbb C$. Then
the natural inclusion $Hom_{\Cal B}(X,\Bbb C)\hookrightarrow
H_{\Cal B,w}(X)^*$ is an isomorphism of vector spaces.
\endproclaim

\prf  As $\Cal B$ is countable, we can assume there exists
a sequence $\{f_k\in\Cal B\}$ satisfying the following properties:
(i) for any $f\in\Cal B$ there exists $k$ with $f\le f_k$, (ii)
$2f_k\le f_{k+1}$ for each $k\ge 1$. Now there is an obvious
inclusion of vector spaces. Suppose that this inclusion is not
also a surjection. Then there exists $\phi\in H_{\Cal B,w}(X)^*$
such that for each $k\ge 1$ there is an $x_k\in X$ with
$|\phi(x_k)|\ge f_k(w(x_k))$. Because $w$ is proper, we may
additionally assume that $w(x_{k+1})
> w(x_k)$ for each $k$. Fix such a choice of $x_k$ for each positive
integer $k$. Let $u_k = \phi(x_k)/|\phi(x_k)|$. Also for each $k$
choose $\alpha_k$ with $\alpha_k/|\alpha_k| = \overline{u_k}$ and
$|\alpha_k| = f_k(w(x_k))^{-1}$. Then
$$
\left|\sum_k\alpha_k x_k\right|_{f_n} \le A_n + \sum_{k
> n}|\alpha_k|\left| x_k\right|_{f_n}\le A_n + \sum_{m\ge 1} 2^{-m} <\infty
$$
for all $n\ge 0$, which by property (i)
implies $\sum_k\alpha_k x_k\in H_{\Cal B,w}(X)$. On the other
hand,
$$
\left|\phi\left(\sum_k \alpha_k x_k\right)\right| = \sum_k
|\alpha_k| |\phi(x_k)|\ge \sum_k 1 =\infty
$$
which is a contradiction. Hence no such $\phi$ can exist.\hfill //
\endpf
\vskip.3in

\subhead 1.3. Topological vs. bornological algebras, and weighted
simplicial sets
\endsubhead
\vskip.1in

Throughout this section we assume $\Cal B$ is countable. There are
apriori different definitions of the topological (bornological)
Hochschild and cyclic homology and cohomology groups associated to
the rapid decay algebras $H_{\Cal B,L}(G)$. We give a brief
description of each, and then show their equivalence. As we have
observed above, $H_{\Cal B,L}(G)$ is a Fr\'echet algebra, hence a
complete, locally convex topological algebra. To do homological
algebra in this setting, we recall some basic definitions.

\exsp \noindent A \underbar{topological $A$-module} is a complete
locally convex space, together with a jointly continuous
$A$-module structure. A topological $A$-module is said to be
\underbar{topologically projective} if it is a direct summand of a
topological $A$-module of the form $A \ttensor E$, for $E$ a
complete locally convex space, and where $\ttensor$ is the
complete projective topological tensor product.

Let $M$ be a topological $A$-module. The proper framework for
defining topological homology and cohomology groups is that of
\underbar{relative} homological algebra, where continuous
surjections of topological $A$-modules are required to be $\Bbb
C$-split [T1, \S 2]. In this setting, a  \underbar{topological
projective resolution of $M$} is an exact sequence of topological
$A$-modules:
$$
0 \leftarrow M \overset{\epsilon}\to{\leftarrow} M_0
\overset{\partial_1}\to{\leftarrow}M_1
\overset{\partial_2}\to{\leftarrow} M_2
\overset{\partial_3}\to{\leftarrow}...
$$
such that each $M_i$ is topologically projective, all of the maps
are continuous $A$-module maps, and the resolution has a
continuous linear homotopy $s_i : M_i \to M_{i+1}$ such
that $\partial_{i+1} s_i + s_{i-1} \partial_i = id_{M_i}$.

\exsp \noindent Given a topological projective resolution $P_{\bullet}$
of $M$ and any topological $A$-module $N$, one can form the complexes
$$
\gathered
P_{\bullet}\underset A\to\tensor N = 0\leftarrow M_0\underset A\to\ttensor N
\overset{\partial_1\otimes Id}\to\leftarrow M_1\underset A\to\ttensor N
\overset{\partial_2\otimes Id}\to{\leftarrow} M_2\underset A\to\ttensor N
\overset{\partial_3\otimes Id}\to{\leftarrow}...\\
Hom_A(P_{\bullet},N) =  0 \rightarrow Hom_A(M_0, N) \overset{\partial_1^*}\to{\rightarrow}
Hom_A(M_1,N) \overset{\partial_2^*}\to{\rightarrow} Hom_A(M_2,N)
\overset{\partial_3^*}\to{\rightarrow} ...
\endgathered
\tag1.3.1
$$
where $Hom_A(M',M'')$ denotes the internal $Hom$-functor evaluated
on $(M',M'')$ (i.e., the vector space of continuous $A$-linear
maps from $M'$ to $M''$). As in the algebraic case, one has
$Tor^A_*(M,N) := H_*(P_{\bullet}\underset A\to\ttensor N)$,
$Ext_A^*(M, N) := H^*(Hom_A(P_{\bullet},N) )$. These definitions
are independent of the particular topological projective
resolution by the usual arguments [Ta]. \vskip.1in

{\bf\underbar{Note 1.3.2}} Given a topological chain complex
$(C_*,d_*)$ with continuous boundary maps, the
\underbar{unreduced} $n^{th}$ homology group is typically meant to
mean $ker(d_n)/im(d_{n+1})$, while the \underbar{reduced} $n^{th}$
homology group is computed as $ker(d_n)/\ov{im(d_{n+1})}$, where
$\ov{im(d_{n+1})}$ denotes the closure of $im(d_{n+1})$ in $C_n$.
The unreduced groups are \lq\lq algebraically correct\rq\rq, while
the reduced groups are \lq\lq topologically correct\rq\rq, in the
sense that the induced subquotient topology is Hausdorff if it was
to begin with on the chain level. In this paper, homology groups
of topological chain complexes will always be unreduced, and the
homology groups viewed simply as vector spaces, without any
inherited topology. The same remarks apply to topological cochain
complexes and cohomology. We will return to this point briefly
later on. 
\vskip.2in

Given a complete, locally convex topological algebra $A$ with unit
and a topological $A$-bimodule $M$, one may define the topological
Hochschild homology of $A$ with coefficients in $M$ as
$HH_*^t(A,M) := Tor^{A\tensor A^{op}}_*(A,M)$. These derived
functors may alternatively be computed as the homology of the
topological Hochschild complex $(C_*^t(A,M),b_*)$ where $C_n(A,M)
:= M\tensor A^{\tensor n}$, with $b_n$ given as above in (1.1).
When $M = A$, the topological cyclic and periodic cyclic complexes
are defined similarly, by taking the algebraic definition in (1.1)
and replacing each occurrence of $(C_*(A,A),b_*)$ with
$(C_*^t(A,A),b_*)$. For cohomology, one defines $HH^*_t(A,M)$ as
$Ext_{A\tensor A^{op}}^*(A,M)$. Again, by convention, the
topological Hochschild cohomology of $A$ - $HH^*_t(A)$ - which is
dual to $HH_*^t(A)$, may be computed as the cohomology of
$Hom(C_*^t(A),\Bbb C) =$ the cocomplex of continuous functionals
on $C_*^t(A)$, which in turn identifies with the cocomplex
$C^*_t(A,A^*)$ ($A^* =$ the continuous dual of $A$; cf. (1.1.1),
[Co1]). Similarly, $HC^*_t(A)$ resp. $HPer^*_t(A)$ is computed as
the cohomology of the cocomplex $Hom(CC_*^t(A),\Bbb C)$ resp.
$Hom(CPer_*^t(A),\Bbb C)$. 

An alternative (and sometimes preferable) setting for doing
homological algebra is the category of bornological algebras and
vector spaces, as shown by Meyer in [Me2]. Briefly, a bornology on
a vector space consists of a collection of bounded subsets. The
bornology is \underbar{convex} if every bounded subset is
contained in an absolutely convex bounded subset, and
\underbar{complete} if every bounded subset is contained in a
complete disk. Bornologies are assumed to be complete and convex
(the original reference for these definitions are [HN1], [HN2],
but [M2, \S 2] will be sufficient for our purposes). As in the
topological case, resolutions in this category are assumed to be $\Bbb C$-split;
moreover the (projective) complete tensor product of two
bornological vector spaces $V$ and $W$ exists and is also denoted
$V\tensor W$ when there is no confusion. Working in the category
of bornological unital algebras and bornological bi-modules over
those algebras, one may proceed as above and define bornological
Hochschild, cyclic and periodic cyclic (co)homology in terms of
the appropriate derived functors. The corresponding groups will be
indicated with a subscript or superscript \lq\lq bor\rq\rq in place
of \lq\lq t\rq\rq. As these functors are originating in different
categories, they are different functors. However, our primary
interest will be the rapid decay algebras defined in (1.1). These
are Fr\'echet algebras, which can be viewed either topologically
or bornologically.

\proclaim{\bf\underbar{Lemma 1.3.3}} For all $(G,L)$ and countable
$\Cal B$, there are natural isomorphisms
$$
\gather
F^t_*(H_{\Cal B,L}(G))\cong F^{bor}_*(H_{\Cal B,L}(G))\\
F_t^*(H_{\Cal B,L}(G))\cong F_{bor}^*(H_{\Cal B,L}(G))
\endgather
$$
for $F = HH, HC$ or $HPer$.
\endproclaim

\prf Consider two weighted sets $(X,w)$ and $(X',w')$. As above,
we may associate to these respectively the Fr\'echet spaces
$H_{\Cal B,w}(X)$ and $H_{\Cal B,w'}(X')$. It is then easy to see that in both
the topological and bornological categories there is an
isomorphism of topological resp. bornological vector spaces
$$
H_{\Cal B,w}(X)\tensor H_{\Cal B,w'}(X')\cong H_{\Cal B,w\times w'}(X\times X')
$$
where $(w\times w')(x,x') = w(x) + w'(x')$. Iterating this
isomorphism yields isomorphisms
$$
H_{\Cal B,w}(X)^{\tensor m}\cong H_{\Cal B,w^m}(X^m)
$$
which apply in either category (where $w^m = w\times
w\times\dots\times w$ $m$ times). Specializing to the case $(X,w)=
(G,L)$, we see that for both projective complete tensor products
there is an isomorphism
$$
H_{\Cal B,L}(G)^{\tensor m}\cong H_{\Cal B,L^m}(G^m)
$$
by which we may identify (up to isomorphism) both the topological
and bornological Hochschild complexes with the complex $(\{H_{\Cal
B,L^{n+1}}(G^{n+1})\}_{n\ge 0},b_*)$, where $b_n$ is the unique
additive extension of the Hochschild boundary map defined on
elements $(g_0,g_1,\dots,g_n)$ by the usual formula (1.1.1). With
this isomorphism established for the Hochschild complex, the
isomorphism for the cyclic and periodic cyclic complexes follows.
It also dualizes to prove the corresponding result in
cohomology.\hfill //
\endpf

In the above situation there is a third, simpler description of
the cochain complex. Again, consider $H_{\Cal B,L}(G)$. By Lemma
1.2.4, there is an isomorphism $H_{\Cal B,L}(G)^*\cong Hom_{\Cal
B}(G,\Bbb C)$. Passing to completed tensor powers and taking
continuous duals results in an identification of both
$C_t^*(H_{\Cal B,L}(G))$ and $C_b^*(H_{\Cal B,L}(G))$ with the
cocomplex
$$
C^*_{\Cal B}(\Bbb C[G]) := \left\{f\in C^n(\Bbb C[G])\ {\dsize |}\ \exists
F\in\Cal B\text{ with } |f(g_0,g_1,\dots,g_n)|\le
F\left(\sum_{i=0}^n L(g_i)\right)\right\} _{n\ge 0}\tag1.3.4
$$
The condition above on the right is simply the statement that $f$
is \underbar{$\Cal B$-bounded} (with respect to the norm on $\Bbb
C$). Thus this cocomplex is the subcocomplex of $C^*(\Bbb C[G])$
consisting of $\Cal B$-bounded cochains; we denote its cohomology
by $HH^*_{\Cal B}(\Bbb C[G])$. We note that the $B_*$ map of
section 1.1 is $\Cal B$-bounded for all $\Cal B$, so one can
define the \underbar{$\Cal B$ cyclic cocomplex} as
$$
CC^*_{\Cal B}(\Bbb C[G]) :=
\{f\in CC^n(\Bbb C[G])\ |\ f\text{ is }\Cal B-\text{bounded}\}_{n\ge 0}
$$
with corresponding cohomology groups denoted $HC^*_{\Cal B}(\Bbb C[G])$.
And similarly with periodic cyclic cohomology.
For any inclusion $\Cal B\hookrightarrow \Cal B'$ of bounding sets
there is an evident transformation of cohomology theories
$F^*_{\Cal B}(\Bbb C[G])\to F^*_{\Cal B'}(\Bbb C[G])$ which preserves all of the standard
structure maps relating the different groups; in particular, it induces a map of
Connes-Gysin sequences. To further understand these cohomology groups,
we introduce some additional terminology.

A \underbar{$\Cal B$-simplicial set}  is a simplicial set
$X\hskip-.03in .$ together with weight functions $w_n:X_n\to\Bbb
R_+$ such that the face and degeneracy maps are $\Cal B$-bounded.
As before, let $(G,L)$ be a discrete group $G$ equipped with a
proper word-length function $L$. A \underbar{$\Cal B$-$G$-set}
consists of a weighted set $(S,w_S)$ together with a specified
left (resp. right) action of $G$ on $S$ such that there exists an
$h\in\Cal B$ with $w_S(gs)$ (resp. $w_S(sg)$) $\le w_S(s) + h(L(g))$ for all $s\in S$ and
$g\in G$ (in other words, the action is uniformly $\Cal
B$-bounded). A $\Cal B$-$G$-map between two $\Cal B$-(left resp.
right) $G$-sets is a $\Cal B$-bounded map which is
$G$-equivariant; it is a $\Cal B$-$G$-isomorphism if it admits a
two-sided inverse which is also a $\Cal B$-$G$-map. Finally, a
\underbar{$\Cal B$-$G$-bi-set} is a weighted set equipped with
both a left and right action of $G$ which is $\Cal B$-bounded in
the above sense, an isomorphism of two such having the obvious
properties. The following are two basic examples relevant to this
paper. \vskip.2in

{\bf\underbar{Example 1.3.5 (i)}} \underbar{The cyclic bar
construction with coefficients} Let $(G,L)$ be as above and let
$(S,w_S)$ be a $\Cal B$-$G$-bi-set. Then $N.^{cy}(G,S) =
\{n\mapsto S\times G^{n}\}_{n\ge 0}$ with
$$
\gather
\partial_0(s,g_1,\dots,g_n) = (sg_1,g_2,\dots,g_n),\\
\partial_i(s,g_1,\dots,g_n) = (s,g_1,\dots,g_ig_{i+1},\dots,g_n),\ 1\le i\le n-1,\\
\partial_n(s,g_1,\dots,g_n) = (g_ns,g_1,g_2,\dots,g_{n-1}),\\
s_j(s,g_1,\dots,g_n) = (s,\dots,g_j,1,g_{j+1},\dots,g_n)
\endgather
$$
The simplicial weight is given by $w_n(s,g_0,\dots,g_n) = w_S(s)
+\sum_{i=1}^n L(g_i)$. The condition on $(S,w_S)$ implies
$N.^{cy}(G,S)$ is a $\Cal B$-simplicial set. When $(S,w_S) =
(G,L)$, we denote this $\Cal B$-simplicial set simply by
$N.^{cy}(G)$ \vskip.1in

{\bf\underbar{Example 1.3.5 (ii)}} \underbar{The bar resolution
with coefficients} Let $(S,w_S)$ be a right $\Cal B$-$G$-set.
Recall that the non-homogeneous bar resolution of $G$ is  $EG. =
\{n\mapsto G^{n+1}\}_{n\ge 0}$ with
$$
\gather
\partial_i[g_0,\dots,g_n] = [g_0,\dots,g_ig_{i+1},\dots,g_n], 0\le i\le n-1,\\
\partial_n[g_0,\dots,g_n] = [g_0,\dots,g_{n-1}],\\
s_j[g_0,\dots,g_n] = [g_0,\dots,g_j,1,g_{j+1},\dots, g_n]
\endgather
$$

The simplicial weight function on $EG.$ is given by
$w([g_0,\dots,g_n]) = \sum_{i=0}^n L(g_i)$. The left $G$-action is
given, as usual, by $g[g_0,g_1,\dots,g_n] = [gg_0,g_1,\dots,g_n]$.
Note that with respect to the given weight function and action of
$G$, $EG.$ is a $\Cal B$-simplicial $G$-set for any $\Cal B$.
Viewing $S$ as a simplicial set with trivial simplicial structure,
we can form the simplicial set $S\times EG.$ with diagonal
simplicial structure. Then $ S\underset G\to\times EG.$ is the
quotient of $S\times EG.$ by the relation
$(sg,[g_0,\dots,g_n])\simeq (s,[gg_0,\dots,g_n])$. This is a $\Cal
B$-simplicial set where the simplicial weight is given by
$$
w([s,[g_0,\dots,g_n]]) :=
\underset{g\in G}\to{inf} \left(w_S(sg^{-1}) + L(gg_0) +  \sum_{i=1}^n L(g_i)\right)
$$
\vskip.2in

To a given $\Cal B$-simplicial set $(X.,w.)$ one can associate the
cochain cocomplex $C^*_{\Cal B}(X,\Bbb C)$ of $\Cal B$-bounded
singular cocains, where
$$
C^n_{\Cal B,w.}(X.,\Bbb C) = \{f:X_n\to\Bbb C\ {\dsize |}\ f\text{
is } \Cal B-\text{bounded with respect to the norm on }\Bbb C\}
$$
This is naturally a subcocomplex of the usual cocomplex of
singular cochains on $X.$ with coefficients in $\Bbb C$. The
\underbar{$\Cal B$-cohomology of $X.$} (with respect to the weight
function $w.$) is $H^*_{\Cal B,w.}(X.) := H^*(C^*_{\Cal
B,w.}(X.,\Bbb C))$. The following is now evident.

\proclaim{\bf\underbar{Lemma 1.3.6}} Under the natural isomorphism
$C^*(\Bbb C[G])\cong C^*(N.^{cy}(G),\Bbb C)$, we have a natural
identification
$$
C_t^*(H_{\Cal B,L}(G))\cong C_b^*(H_{\Cal B,L}(G))\cong C^*_{\Cal
B,w.}(N.^{cy}(G),\Bbb C)
$$
for all $(G,L)$ and countable $\Cal B$.
\endproclaim

For convenience, we will denote the groups $C^*_{\Cal
B,w.}(N.^{cy}(G),\Bbb C)$ simply by $C^*_{\Cal
B}(N.^{cy}(G))$.
\vskip.3in

\subhead 1.4. The decomposition
\endsubhead
\vskip.1in

As in the case of the group algebra, the Hochschild , cyclic and
periodic cyclic (co)homologies of the afore-mentioned rapid decay
algebras admit a decomposition indexed on $<G> =$ the set  of
conjugacy classes of $G$. The starting point is the following
simple

\proclaim{\bf\underbar{Proposition 1.4.1}} For any class of
bounding functions $\Cal B$, there is an isomorphism of $\Cal
B$-simplicial sets
$$
N.^{cy}(G)\cong \underset{x\in <G>}\to\coprod (N.^{cy}(G)_x)\cong
\underset{x\in <G>}\to\coprod S_x\underset{G}\to\times EG.
$$
(cf. (1.3.5 (i), (ii) above), where $(N.^{cy}(G)_x)$ denotes the
$\Cal B$-subsimplicial set of $N.^{cy}(G)$ with $(N.^{cy}(G)_x)_n =
\{ (g_0,\dots,g_n)\ |\ <\prod_{i=0}^n g_i> = x\}$, and $S_x =
\{g\in G\ |\ <g> = x\}$ with weight function induced by $L$ and
the inclusion into $G$. The right $G$-action on $S_x$ is given by
$h(g) := g^{-1}h g$ (conjugation by $g^{-1}$).
\endproclaim

\prf This equivalence is well-known, and the first decomposition
is clearly a $\Cal B$-simplicial isomorphism. It remains to note
that the maps
$$
\gather
(N.^{cy}(G)_x)_n\ni (g_0,g_1,\dots,g_n)\mapsto
[(\prod_{i=0}^n g_i),[g_0,\dots,g_n]] =
[(\prod_{i=0}^n g_i)^{g_0^{-1}},[1,\dots,g_n]]\in S_x \underset{G}\to\times EG_n\\
S_x \underset{G}\to\times EG_n\ni [g,[1,\dots,g_n]]\mapsto
((\prod_{i=1}^ng_i)^{-1}g,g_1,\dots,g_n)
\endgather
$$
are $\Cal B$-bounded for any $\Cal B$.\hfill //
\endpf

Together with the results of the previous section, we have the
following immediate consequence.

\proclaim{\bf\underbar{Corollary 1.4.2}} For any conjugacy class
$x$ the decomposition of (1.4.1) induces factorizations in homology and cohomology
$$
\gather
\underset{x\in <G>}\to{{\oplus}}F_*^t(H_{\Cal B,L}(G))_x\inj F_*^t(H_{\Cal B,L}(G))\to
\underset{x\in <G>}\to{{\prod}}F_*^t(H_{\Cal B,L}(G))_x \\
\underset{x\in <G>}\to{{\oplus}}F^*_t(H_{\Cal B,L}(G))_x\inj F^*_t(H_{\Cal B,L}(G))\to
\underset{x\in <G>}\to{{\prod}}F^*_t(H_{\Cal B,L}(G))_x 
\endgather
$$
for $F = HH, HC, HPer$, and with the maps in the Connes-Gysin
sequence preserving this decomposition. These summands are preserved by homomorphisms coming from an inclusion of algebras $H_{\Cal B,L}(G)\hookrightarrow H_{\Cal B',L}(G)$ (induced by the relation 
$\Cal B'\prec \Cal B$). Moreover there exists, for
each $x$, isomorphisms
$$
HH^*_t(H_{\Cal B,L}(G))_x\cong H^*_{\Cal B}(
S_x\underset{G}\to\times EG.)
$$
\endproclaim

We note that the the homomorphism $F_*^t(H_{\Cal B,L}(G))\to
\underset{x\in <G>}\to{{\prod}}F_*^t(H_{\Cal B,L}(G))_x $ may \underbar{not} be injective (with the same remark holding for the cohomology groups). In this paper, however, we will only be concerned
with the summands indicated in the sums on the left.

Fix a conjugacy class $x$ and group element $h$ with $<h> = x$,
let $G_h $ denote the centralizer of $h$ in $G$, and
$G_h\backslash G$ the corresponding right coset space. This is a
weighted set, with weight $\ov{L}(G_hg) := \underset{g'\in
G_h}\to{min} L(g'g)$. The right $G$-action on $G_h\backslash G$ is
the obvious one.

\proclaim{\bf\underbar{Proposition 1.4.3}} The natural map
$\pi_x:G_h\backslash G\to S_x$ given by $G_hg\mapsto g^{-1}hg$ is
a $G$-equivariant isomorphism of right $G$-sets, and $\Cal
B$-bounded.
\endproclaim

\prf It is obviously a $G$-equivariant isomorphism. It is also
linearly bounded in terms of the respective weight
functions.\hfill//
\endpf

Unfortunately, and here is the key point, this map is
\underbar{not} a $\Cal B$-$G$-isomorphism of right $\Cal B$-$G$-sets
in general. In order to proceed beyond Cor. 1.4.2 and get
an identification analogous to Theorem 1.1.2 (ii), we will need to
impose conditions on $(G,L)$. We say that
$(G,L)$ has a \underbar{$\Cal B$-solvable conjugacy-bound} at the conjugacy class  $x = <h>\in <G>$ if there exists $f\in\Cal B$ with the
property: for each $y\in S_x$, there exists $g = g_y\in G$ with $y
= h^{g_y}$ and $L(g_y)\le f(L(y))$. In other words, there exists a
map $\psi_x:S_x\ni y\mapsto [g_y]\in G_h\backslash G$ which is an
inverse to $\pi_x$ and is $\Cal B$-bounded. As $\pi_x$ is
$G$-equivariant, so is $\psi_x$. The pair $(G,L)$ has a $\Cal B$-solvable conjugacy-bound if it does so at each non-elliptic conjugacy class $x$.

\proclaim{\bf\underbar{Lemma 1.4.4}} The pair
$(G,L)$ has a $\Cal B$-solvable conjugacy-bound at $x$ iff $\pi_x$ is a $\Cal
B$-$G$-isomorphism of right $\Cal B$-$G$-sets. In this case, there
is a natural isomorphism of $\Cal B$-simplicial sets
$$
(\pi_x).: G_h\backslash G \underset{G}\to\times
EG.\overset\cong\to{\underset\Cal B\to\longrightarrow} S_x
\underset{G}\to\times EG.
$$
which induces an isomorphism on $\Cal B$-bounded cohomology.
\endproclaim

\prf Given the definitions of the respective simplicial weight
functions, it suffices to check that $(\pi_x)_1$ is a $\Cal
B$-isomorphism, which is implied by the fact $\pi_x$ is a $\Cal
B$-isomorphism.\hfill //
\endpf

\proclaim{\bf\underbar{Proposition 1.4.5}} For each $h$ the
inclusion $G_h\hookrightarrow G$ yields an inclusion of $\Cal
B$-simplicial sets

$$
BG_h = (G_h\backslash G) \underset{G_h}\to\times
E(G_h).\hookrightarrow G_h\backslash G \underset{G}\to\times EG.
$$
which induces an isomorphism on $\Cal B$-bounded cohomology
groups.
\endproclaim

\prf This result is shown in [Og2] for $\Cal B = \Cal P$, but the
same argument works in this more general case (this point is also
noted in [Me1] for $\Cal P$ and $\Cal E$).\hfill //
\endpf

Combining 1.4.2 and 1.4.5 we have

\proclaim{\bf\underbar{Corollary 1.4.6}} If the pair
$(G,L)$ has a $\Cal B$-solvable conjugacy-bound, then for each non-elliptic conjugacy class $x
= <h>$, there is an isomorphism of cohomology groups
$$
HH^*_t(H_{\Cal B,L}(G))_x\cong H^*_{\Cal B}(
S_x\underset{G}\to\times EG.)\cong H^*_{\Cal B}(BG_h;\Bbb C)
$$
\endproclaim

We observe that having a $\Cal B$-solvable conjugacy-bound is a local condition, and strictly
weaker than requiring that the $G$-map
$$
\diagram
\underset{x=<h>\in <G>}\to\coprod (G_h\backslash G)
\xto[0,3]^{\underset{x=<h>\in <G>}\to\coprod \pi_x}
& & &\underset{x=<h>\in <G>}\to\coprod (S_x)
\enddiagram
$$
is a $\Cal B$-$G$-isomorphism. Given $G$, it is also possible that
certain conjugacy classes have a $\Cal B$-solvable conjugacy-bound, while others are
not. To illustrate, $<1>$ has a $\Cal B$-solvable conjugacy-bound for any $G$ and
bounding class $\Cal B$, for trivial reasons. This statement is
false in general for conjugacy classes $x\ne <1>$.

Borrowing some terminology from [Me2], we say that $(G,L)$ is
\underbar{$\Cal B$-isocohomological} if the transformation
$H^*_{\Cal B}(BG;\Bbb C)\to H^*(BG;\Bbb C)$ (induced by the
inclusion $\Cal B\hookrightarrow \Cal B_{max}$) is an isomorphism.
Putting this all together, we have the following result.

\proclaim{\bf\underbar{Theorem 1.4.7}} Let $(G,L)$ be a discrete
group equipped with proper word-length function $L$. Suppose there
exists a countable bounding class $\Cal B$ such that
\vskip.05in

(i) Each non-elliptic conjugacy class $x$ has a $\Cal B$-solvable conjugacy-bound, and
\vskip.05in
(ii) For each non-elliptic $x = <h>\in <G>$, $G_h$ is $\Cal
B$-isocohomological.
\vskip.05in

Then for each such $x$ there are isomorphisms
$$
\gather
HH^*_t(H_{\Cal B,L}(G))_x = HH^*_{\Cal B}(\Bbb C[G])_x\cong
HH^*(\Bbb C[G])_x = H^*(BG_h;\Bbb C)\\
HC^*_t(H_{\Cal B,L}(G))_x = HC^*_{\Cal B}(\Bbb C[G])_x\cong
HC^*(\Bbb C[G])_x
\endgather
$$
for each $x=<h>\in<G>$.
\endproclaim

\prf As we have seen, $HH^*_{\Cal B}(\Bbb C[G])_x\cong H^*_{\Cal
B}(BG_h;\Bbb C)$ when $x$ has a $\Cal B$-solvable conjugacy-bound. Condition (ii)
then implies $H^*_{\Cal B}(BG_h;\Bbb C)\cong H^*(BG_h;\Bbb C)$.
The statement for cyclic cohomology follows by a comparison
of Connes-Gysin sequences.\hfill //
\endpf

As with non-bounded Hochschild and cyclic cohomology groups, it is
relatively easy to show that when $x$ has a $\Cal B$-solvable conjugacy-bound and is
\underbar{elliptic}, there is an isomorphism
 $$
 HC^*_{\Cal B}(\Bbb C[G])_x\cong H^*_{\Cal B}(BG_h;\Bbb C)\otimes HC^*(\Bbb C)
 $$
where the $S^*$-map on the right is given by $Id\otimes S^*_{\Bbb
C}$. For $x$ non-elliptic, the situation is more delicate. The
next lemma follows by a straightforward adoption of the method of
Nistor [Ni1] (cf. [Ji2]). As it is not needed for the rest of the
paper, we state the it without proof.

\proclaim{\bf\underbar{Lemma 1.4.8}} Given $(G,L)$ satisfying the
conditions of Theorem 1.4.7, and $x = <h>\in <G>$ a non-elliptic
class. Let $C_h\subset G_h$ denote the infinite cyclic subgroup
generated by $h$, and let $L_h = L|_{C_h}$. Let $\Bbb Z$ denote
the infinite (multiplicative) cyclic group generated by $t$, with
$L_{\Bbb Z}(t^n) = n$. If the isomorphism $(\Bbb Z,L_{\Bbb Z})\to
(C_h,L_h);\ t\mapsto h$ is a $\Cal B$-isomorphism, then
$HC^*_t(H_{\Cal B,L}(G))_x\cong H^*_{\Cal B}(N_h;\Bbb C)$, where
$N_h = G_h/C_h$ equipped with the quotient word-length function
coming from $G_h$.
\endproclaim

Our main application of Theorem 1.4.7 is

\proclaim{\bf\underbar{Corollary 1.4.9}} Suppose $G$ is in the
class $\cC(\Bbb C)$, and $(G,L)$ satisfies the conditions of
Theorem 1.4.7. Then $(G,L)$ is in the class $\cC_{\Cal
B}(\Bbb C)$. Consequently the conjecture $\Cal B$-SrBC holds for $G$.
\endproclaim

 In the next section we exhibit certain classes of groups to which this Corollary applies.
\vskip.4in

\head 2. Geometric groups satisfying $\Cal B$-SrBC
\endhead
\vskip.2in

\subhead 2.1. Nilpotent and hyperbolic groups
\endsubhead
\vskip.1in

In previous work, the first author has shown
$(G,L)$ satisfies the $\Cal P$-nilpotency condition for
any proper word-length function $L$ on $G$ when $G$ is of polynomial growth [Ji1] or
word-hyperbolic [Ji3]. We show that for word hyperbolic groups or finitely-generated
torsion free two-step nilpotent groups, these results can be efficiently recovered
(and extended) using the machinery of section 1.
\vskip.1in

\proclaim{\bf\underbar{Proposition 2.1.1}} Let $G$  be a finitely-generated torsion free
two-step nilpotent
group, and $L$ a proper word-length function on $G$. Then for all bounding families
$\Cal B\supseteq\Cal P$
\vskip.05in

(i) Each conjugacy class $x$ has a $\Cal B$-solvable conjugacy-bound,
\vskip.05in

(ii) The centralizer $G_h$ for $h \in x$ is $\Cal B$-isocohomogical for
each conjugacy class $x$, and
\vskip.05in

(iii) For each non-elliptic $x = <h>$, $G_x/(h)$ has finite rational
cohomological dimension ($r.c.d.< \infty$).
\vskip.05in

Hence $G$ satisfies $\Cal B$-nilpotency for all $\Cal B\supseteq \Cal P$.
\endproclaim

\prf For (i), we note that polynomial solvability of the conjugacy problem
implies the existence of a $\Cal B$-solvable conjugacy-bound for all $\Cal B\supseteq \Cal P$. 
Note that by [Bl] the conjugacy problem for a finitely-generated nilpotent
group is solvable.  We must then show that there is a polynomial solution.
Consider the torsion-free nilpotent group $G$ given by the central
extension $C \rightarrow G \rightarrow N$, with $C$ and $N$ free abelian
groups, with bases $\{e'_1, \ldots, e_n\}$ and $\{ f'_1, \ldots, f'_m\}$ respectively.
Let $g$ and $h$ be conjugate elements in $G$, with
$ug = hu$.  We can write these elements as 
$g = f_1^{x_1} f_2^{x_2} \ldots f_m^{x_m} c$, 
$h=f_1^{y_1} f_2^{y_2} \ldots f_m^{y_m} c'$, and 
$u=f_1^{z_1} f_2^{z_2} \ldots f_m^{z_m}$ with $x_i$, $y_i$, and $z_i$ integers,
$c$ and $c'$ elements of $C$. We can also write
$c=e_1^{a_1} \ldots e_n^{a_n}$, and $c'=e_1^{b_1} \ldots e_n^{b_n}$, for integers $a_i$ and $b_j$.
( Here we have identified $C$ with its image in $G$,
and $f_i$ is the image of $f'_i$ in $G$ under a fixed cross-section $N \rightarrow G$.
$f'_i f'_j = f'_j f'_i$ so that $f_i f_j = f_j f_i C_{i,j}$ for some $C_{i,j} \in C$.
From $ug = hu$ we have $f_1^{z_1} f_2^{z_2} \ldots f_m^{z_m}f_1^{x_1} f_2^{x_2} \ldots f_m^{x_m} c = 
f_1^{y_1} f_2^{y_2} \ldots f_m^{y_m} f_1^{z_1} f_2^{z_2} \ldots f_m^{z_m} c'$.  Applying
the conjugation properties we obtain
$f_1^{x_1+z_1} f_2^{x_2+z_2} \ldots f_m^{x_m+z_m}C_{m,1}^{x_m z_1} \ldots C_{1,1}^{x_2 z_1}C_{m,2}^{x_m z_2} \ldots
C_{m,m-1}^{x_m z_{m-1}}c = 
f_1^{z_1 + y_1} f_2^{z_2+y_2} \ldots f_m^{z_m+y_m}C_{m,1}^{z_m y_1} \ldots C_{1,1}^{z_2 y_1}C_{m,2}^{z_m y_2} \ldots
C_{m,m-1}^{z_m y_{m-1}} c'$  Moving all of the $f_i$'s to tone side we obtain
$f_1^{x_1-y_1} f_2^{x_2-y_2} \ldots f_m^{x_m-y_m}C_{m,1}^{x_m z_1-y_1 z_m} \ldots C_{1,1}^{x_2 z_1-y_1 z_2}C_{m,2}^{x_m z_2-y_2 z_m} \ldots C_{m,m-1}^{x_m z_{m-1}-y_{m-1} z_m} c {c'}^{-1} = e$
As $\{f'_1, \ldots, f'_m\}$ are free abelian generators for $N$, we see that $x_i = y_i$ for all $i$.
The last equation reduces to 
$C_{m,1}^{x_m z_1-x_1 z_m} \ldots C_{1,1}^{x_2 z_1-x_1 z_2}C_{m,2}^{x_m z_2-x_2 z_m} \ldots C_{m,m-1}^{x_m z_{m-1}-x_{m-1} z_m} c {c'}^{-1} = e$.  From this we obtain a linear system to determine the $z_i$ in terms 
of the $x_j$s, $a_j$s, and $b_j$s.  As $g$ and $h$ are conjugate, this system has a solution.
Thus there is a solution which is bounded by a polynomial of the absolute values of the entries in the
coefficient matrix of this system.  Moreover as $C$, $N$, and $G$ are fixed, this matrix, and hence the
polynomial itself, is independent of $g$ and $h$.

As $G$ is nilpotent, so is $G_x$ for all $x$. 
Choose an element
$z \in C$ with $ord(z) = \infty$. Let $C_z$ denote the
infinite cyclic subgroup generated by $z$. Associated to the
short-exact sequence of groups $\Bbb Z = C_z\inj G\surj \ov{G}$ is
a Serre spectral sequence in $\Cal B$-cohomology (the case $\Cal B
= \Cal P$ is done in detail in [Og2]. The more general case follows
by similar reasoning, with similar restrictions - cf. [Og3]). 
By a spectral sequence comparison together with induction
on $n$, we see that nilpotent groups are $\Cal B$-isocohomological for all 
$\Cal B\supseteq \Cal P$. Finally, any
f.g. nilpotent group $G'$ satisfies $r.c.d.(G')< \infty$; as
subquotients of f.g. nilpotent groups are again f.g. nilpotent,
this implies (iii).\hfill //
\endpf


\proclaim{\bf\underbar{Proposition 2.1.2}} Let $G$  be a
word-hyperbolic group, and $L$ a proper word-length function on
$G$. Then for all bounding families $\Cal B\supseteq\Cal L$
\vskip.05in

(i) Each conjugacy class $x$ has a $\Cal B$-solvable conjugacy-bound,
\vskip.05in

(ii) The centralizer $G_h$ is $\Cal B$-isocohomogical for
each element $h$ in the non-elliptic conjugacy class $x$, and
\vskip.05in

(iii) For each non-elliptic class $x= <h>$, $G_h/(h)$ has $r.c.d.< \infty$.
\vskip.05in

Hence $G$ satisfies $\Cal B$-nilpotency for all $\Cal B\supseteq
\Cal L$.
\endproclaim

\prf In the case of word-hyperbolic groups, the situation is  even
nicer. 
If $x$ is a non-elliptic conjugacy
class, then the conjugacy-bound in this case is linear [Gr], [Bu]. Moreover, $G_h$ is virtually cyclic and isometrically embedded in $G$ for any $h \in x$, implying it is $\Cal
B$-isocohomological [Ji1], [Og2]. It also implies the quotient $G_{h}/(h)$ is finite, so $r.c.d. =
0$.\hfill //
\endpf

\vskip.3in

\subhead 2.2. Relatively hyperbolic groups
\endsubhead
\vskip.1in

We begin by recalling the basic setup for relatively hyperbolic
groups. Thus, let $G$ be a discrete group, and
$\{H_\lambda\}_{\lambda \in \Lambda}$ a family of subgroups of
$G$. Set $\RH= \bigcup_{\lambda \in \Lambda} H_\lambda \setminus
\{ 1_G \}$.  A subset $X \subset G$ is a relative generating set
of $G$, with respect to $\{H_\lambda\}_{\lambda \in \Lambda}$, if
$G = < X \cup \RH >$. In this case, $G$ can be considered
as the quotient of a free product $F = F(X) * \left( *_{\lambda
\in \Lambda} H_\lambda \right)$, where $F(X)$ is the free group
with basis $X$.

Let $\rel \subset F$ be such that the normal closure of $\rel$ is
the kernel of the quotient map. In this case $G$ has a relative
presentation with respect to $\{H_\lambda\}_{\lambda \in \Lambda}$
given by $< X, \RH | R = 1, R \in \rel>$.  $G$ is finitely
relatively presented if $\rel$ is finite.

Let $W$ be a path in $F$ representing $1_G \in G$.  Then $W$ can
be represented in the following form: $W = \Pi_{i=1}^{k} f_i
R_i^{\pm 1}f_i^{-1}$, where the $f_i \in F$ and $R_i \in \rel$. We
can then define the relative area of such a path: $area^{rel}(W) =
\min \{ k \in \Cal N | W \text{ can be written as } \Pi_{i=1}^{k}
f_i R_i^{\pm 1}f_i^{-1} \text{ where } f_i \in F \text{ and } R_i
\in \rel \}$.

$G$ is \underbar{relatively hyperbolic} with respect to
$\{H_\lambda\}_{\lambda \in \Lambda}$ if $G$ is finitely
relatively presented with respect to $\{ H_\lambda\}_{\lambda \in
\Lambda}$ and there is a constant $C$ such that for any word $W
\in F$ representing $1_G \in G$ we have that $area^{rel}(W) \leq C
\| W \|$, where $\|W\|$ is the length of $W$ in $F$ \cite{Os}.  In
particular if $G$ is relatively hyperbolic with respect to $\{
H_\lambda \}_{\lambda \in \Lambda}$, then the Cayley graph
$\Gamma( G, X \cup \RH)$ is a hyperbolic metric space, in the
graph metric.

An alternate description is given by Farb in [Fa].  Let $G$ be a finitely
generated group, and let $\{H_\lambda\}_{\lambda \in \Lambda}$ be a 
family of finitely generated subgroups of $G$. Fix $X$, a finite symmetric 
set of generators of $G$. To the corresponding Cayley graph $\Gamma(G,X)$ 
add a vertex $v(gH_\lambda)$ for each left coset of each $H_\lambda$ in the 
family, and connect $v(gH_\lambda)$ with the vertex corresponding to each
element $x \in gH_\lambda$ by an edge of length $1/2$.  The resulting graph,
$\hat \Gamma(G,X; \{H_\lambda\}_{\lambda \in \Lambda})$ is called the 
coned-off Cayley graph of $G$ relative to $\{H_\lambda\}_{\lambda \in \Lambda}$.  
In Farb's notation, the group $G$ is weakly relative hyperbolic with respect 
to $\{H_\lambda\}_{\lambda \in \Lambda}$ if 
$\hat \Gamma(G,X; \{H_\lambda\}_{\lambda \in \Lambda})$ is a hyperbolic 
metric space.  As $\hat \Gamma(G,X; \{H_\lambda\}_{\lambda \in \Lambda})$ 
is quasi-isometric to $\Gamma( G, X \cup \RH)$, this property is weaker 
than that defined above.  However, when augmented with the 
bounded coset penetration property ( defined below ), this notion is 
equivalent to relative hyperbolicity.

\proclaim{\bf\underbar{Definition}} 
Let a group $G$ be weakly relatively hyperbolic to a family 
$\{H_\lambda\}_{\lambda \in \Lambda}$ of finitely generated subgroups.  
$(G, \{H_\lambda\}_{\lambda \in \Lambda})$ is said to satisfy the
\underbar{Bounded Coset Penetration (BCP) property} if for each
constant $k$ there exists a constant $c = c(k)$ such that for every
pair of $k$-quasi-geodesics $p$ and $q$ in the coned-off Cayley graph with 
the same endpoints and without backtracking, satisfy\vskip.05in

(i) If $p$ penetrates a coset $gH_\lambda$ and $q$ does not penetrate 
		$gH_\lambda$, then the point at which $p$ enters the coset is
		at most a $\Gamma(G,X)$-distance of $c(k)$ from the point at which $p$ leaves
		the coset;\vskip.05in

(ii) If $p$ and $q$ both penetrate a coset $gH_\lambda$, then the points
		at which $p$ and $q$ enter $gH_\lambda$ are at most a $\Gamma(G,X)$-distance
		of $c(k)$ from each other.  Similarly the points at which $p$ and $q$ exit
		the coset are within a $\Gamma(G,X)$-distance of $c(k)$ from one-another. 
\endproclaim

\proclaim{\bf\underbar{Definition}} 
Let $G$ be a finitely generated group which is relatively hyperbolic
with respect to the collection of subgroups $\{ H_\lambda \}_{\lambda
\in \Lambda}$.  An element of $G$ is \underbar{parabolic} if it is
conjugate to an element in one of the $H_\lambda$.  Otherwise the element
is said to be \underbar{hyperbolic}.
\endproclaim

\proclaim{\bf\underbar{Proposition 2.2.1}}
Let $G$ be a finitely
generated group, $\{ H_\lambda\}_{\lambda \in \Lambda}$ a
collection of nontrivial subgroups of $G$. Suppose that $G$ is
finitely presented with respect to $\{ H_\lambda\}_{\lambda \in
\Lambda}$.  Then the collection $\{ H_\lambda\}_{\lambda \in
\Lambda}$ is finite. In particular, if $G$ is relatively
hyperbolic with respect to a collection $\{ H_\lambda\}_{\lambda
\in \Lambda}$ of nontrivial subgroups, then $card \, \Lambda <
\infty$.
\endproclaim
\prf
This is Corollary 2.48 of [Os].
\hfill //
\endpf

\proclaim{\bf\underbar{Proposition 2.2.2}} Let $G$ be a finitely
generated group with generating set $X$, such that $G$ is
relatively hyperbolic with respect to $\{ H_\lambda \}_{\lambda
\in \Lambda}$, and let $h \in G$ be a nontorsion element which
does not conjugate into any of the $H_\lambda$ (such an element is
said to be hyperbolic). Then the centralizer $G_h$ is finitely
generated, virtually cyclic, and the embedding $G_h \to G$ is a
quasi-isometric embedding.
\endproclaim

\prf
By Theorem 4.19 of [Os], $G_h$ is a strongly relatively 
quasi-convex subgroup of $G$.  Theorem 4.13 of [Os] shows that 
$G_h$ is finitely generated and that the inclusion $G_h \to G$ is a
quasi-isometric embedding.  Theorem 4.16 of [Os] gives that $G_h$ is
a word-hyperbolic group.  As the infinite cyclic subgroup generated 
by $h$ lies in the center of the hyperbolic group $G_h$, we have that
$G_h$ is virtually cyclic.
\hfill //
\endpf

\proclaim{\bf\underbar{Proposition 2.2.3}} Let $G$ be a finitely
generated group with generating set $X$, such that $G$ is
relatively hyperbolic with respect to $\{ H_\lambda \}_{\lambda
\in \Lambda}$, and let $h \in G$ be a nontorsion parabolic element such that
$h \in xH_\lambda x^{-1}$. Then $G_h = x H_{\lambda, x^{-1}hx} x^{-1}$, where
$H_{\lambda, x^{-1}hx}$ is the centralizer of $x^{-1} h x$ in
$H_\lambda$.
\endproclaim

\prf
This is due to the fact that $ H_\lambda $ and $xH_\lambda x^{-1}$ are conjugated by $x$.
\hfill //
\endpf

\proclaim{\bf\underbar{Proposition 2.2.4}} Let $G$ be a finitely
generated group with generating set $X$, such that $G$ is
relatively hyperbolic with respect to $\{ H_\lambda \}_{\lambda
\in \Lambda}$, and let 
$\RH= \bigcup_{\lambda \in \Lambda} H_\lambda \setminus \{ 1_G \}$.
Then for every $k$ there is an
$\epsilon = \epsilon(k)$ such that given any two $k$-quasi-geodesics
without backtracking in $\Gamma(G, X \cup \RH)$, say $p$ and $q$, 
the following conditions hold:\vskip.05in
(1) If $p$ penetrates a coset $gH_\lambda$ and $q$ does not penetrate 
		$gH_\lambda$, then the point at which $p$ enters the coset is
		at most a $\Gamma(G,X)$-distance of $\epsilon(k)$ from the point at which $p$ leaves
		the coset;\vskip.05in
(2) If $p$ and $q$ both penetrate a coset $gH_\lambda$, then the points
		at which $p$ and $q$ enter $gH_\lambda$ are at most a $\Gamma(G,X)$-distance
		of $\epsilon(k)$ from each other.  Similarly the points at which $p$ and $q$ exit
		the coset are within a $\Gamma(G,X)$-distance of $\epsilon(k)$ from one-another.
\endproclaim
Thus $\epsilon$ works for the constant $c$ in the definition of the BCP property.
\prf
This is Theorem 3.23 of [Os].
\hfill //
\endpf

\proclaim{\bf\underbar{Lemma 2.2.5 }}
Let $(H,d)$ be a $\delta$-hyperbolic geodesic metric space.  For
every $k$ there is a constant $\tilde N= \tilde N(k)$ such that if $p$
is a $k$-quasigeodesic with endpoints $x$ and $y$, then $p$
lies within the $\tilde N$-neighborhood of every geodesic $[x,y]$
joining $x$ to $y$ in $H$.  Moreover $\tilde N$ can be bounded by
a polynomial in $k$.
\endproclaim

\prf
This is theorem 1.7 of Part III.H of [BH].
[BH] actually claims only the first statement above, but the 
existance of the polynomial follows from the proof.
In particular, 
$$
\tilde N(k) \le \left( \delta \log_2 \left( 2 k^3 + 6k^2 + 3k + 2 \right) + \delta \log_2 \left[ \delta \log_2 \left( 2 k^3 + 6k^2 + 3k + 2 \right) \right] + 1 \right) \left( k^2 + 1 \right) + \frac{1}{2} \left( 2k^3 + 3k \right)
$$
which is clearly bounded by a polynomial in $k$.\hfill //
\endpf

\proclaim{\bf\underbar{Corollary 2.2.6 }}
Let $(H,d)$ be a $\delta$-hyperbolic geodesic metric space.  For
every $k$ and $R$ there is a constant $N=N(k,R)$ such that if $p$ and $q$
are $k$-quasigeodesics whose starting points are within $R$ from each other,
and whose ending points are within $R$ from each other, then $p$ and $q$
lie within the $N$-neighborhood of one-another.  Moreover $N(k,R)$ can be bounded by
a polynomial in $k$ and $R$.
\endproclaim
\prf
Denote the starting and ending points of $p$ ( resp. $q$ ) by $x$ and $y$
( resp. $x'$ and $y'$ ).  By Lemma 2.2.5, $p$ and $[x,y]$ lie within
an $\tilde N(k)$-neighborhood of one-another.  Similarly for $q$ and $[x',y']$.
Consider the geodesic quadrilateral $x,y,y',x'$.  The sides $[x,x']$ and 
$[y,y']$ each have length less than $R$ by hypothesis.  For any point on $[x,y]$
which is at least $R+2\delta$ from either endpoint, hyperbolicity gives a point
on $[x',y']$ which is within a distance of $2\delta$.  Similarly for points
on $[x',y']$ at least $R + 2\delta$ from either endpoint.  As every point of
$[x,y]$ is within $R+2\delta$ from such a point, it follows that $[x,y]$ and
$[x',y']$ lie in $R + 4\delta$-neighborhoods of one-another.  It follows that
$N(k,R) = \tilde N(k) + R + 2\delta$ works.  Moreover, as $\tilde N(k)$ is
bounded by a polynomial in $k$, we have $N(k,R)$ is bounded by a polynomial
in $k$ and $R$.\hfill //
\endpf

\proclaim{\bf\underbar{Proposition 2.2.7}} 
Let $G$, $\RH$ and $\epsilon$ be as in Proposition 2.2.4.  
Then $\epsilon(k)$ is bounded by a polynomial in $k$.
\endproclaim
\prf
In [Os] the proof of theorem 3.23 bounds $\epsilon(k)$ by the maximum of
three terms, each of which can be expressed as a polynomial in 
$N=N(k)$.  By the above lemma, this is bounded by a polynomial in $k$.

In particular, $\epsilon(k) = \max\{ \epsilon'(k), C'(k), D(k) \}$ where 
$C'(k) = LM \left( 1 + k( 2\epsilon'(k) + 1) + 2 \epsilon'(k) \right)$,
$D(k) = LM \left( 1 + k( 2\epsilon'(k) + 1) + 2 \epsilon'(k) + k + 1 \right)$,
$\epsilon'(k) = 2K(k)^2 L M \left( 4k+1 \right)$, for $L$ and $M$ positive
constants depending only on the finite relative presentation
$< X, \RH | R = 1, R \in \rel>$.  Let $K_0 = N(k, 0)$, where $N(k,R)$
is the constant from Corollary 2.2.6.
Then $K(k) = N(k,K_0) + \frac{1}{2}$ is the $K(k)$ showing up in
the definition of $\epsilon'$.  
It is clear that $\epsilon(k)$ is bounded by a polynomial in $k$.\hfill //
\endpf

The following is, at heart, the Pseudo-Anosov case of Theorem 4.1 of [Ma].
\proclaim{\bf\underbar{Lemma 2.2.8}}
Let $u$ and $v$ be conjugate nontorsion hyperbolic elements of $G$.  Then there exists
a constant $K_h$ and a $g \in G$ with $u = g^{-1} v g$, and 
$\ell_{\hat \Gamma}(\hat g) \le K_h ( \ell_{\Gamma}(u)+ \ell_{\Gamma}(v) )$,
where $\hat g$ is the ``shortened path'' in the relative graph $\hat \Gamma$
corresponding to the element $g \in G$.  Moreover, $K_h$ is independant of 
the choice of $u$ and $v$.
\endproclaim
\prf
By the work in [Os] on translation length of hyperbolic elements, there
exists a $d > 0$ such that for any hyperbolic element $h \in G$ one has
$\ell_{\hat \Gamma}( \hat{h^n} ) \ge |n| d$.
Moreover, for any $x \in G$ we have $d_{\hat \Gamma}(x, h^n x) \ge |n|d$.

Let $\alpha$ be a quasi-axis for the hyperbolic element $h$, acting on $\hat \Gamma$
as a hyperbolic isometry.  That is, $\alpha$ is a bi-infinite geodesic such that
for every $n$, $h^n \alpha$ lies in a $2 \delta$ neighborhood of $\alpha$, where $\delta$
is the hyperbolicity constant of $\hat \Gamma$.

If $x \in G$ then there exists constant $K$ such that 
$d_{\hat \Gamma}(x, \alpha) \le K d_{\hat \Gamma}(x, hx)$, for $K$ independant of $x$ and $h$.
This is Lemma 4.3 of [Ma], which relies only on $\hat \Gamma$ being a hyperbolic metric space,
on which $h$ acts as a hyperbolic isometry.

Thus if $u$ and $v$ are two conjugate hyperbolic elements of $G$, let $g$ be some conjugator
with $u = g^{-1} v g$.  Let $\alpha$ and $\beta$ be the axes of 
$u$ and $v$, respectively, as hyperbolic isometries.  Let $a$ be the closest point
on $\alpha$ to the identity element, $e$, and let $b$ be the closest point on
$\beta$ to $e$.  Applying our above inequality yields 
$d_{\hat \Gamma}(e,a) \le K d_(e,u) = K \ell_{\hat \Gamma}(\hat u)$ and
$d_{\hat \Gamma}(e,b) \le K \ell_{\hat \Gamma}(\hat v)$.  The conjugator $g$
takes the axis $\beta$ to the axis $\alpha$, so $g b \in \alpha$.  The translation
length of $u$ is at most $\ell_{\hat \Gamma}(\hat u)$, so there is an $n$ such that
$d_{\hat \Gamma}(a, u^n g b) \le \ell_{\hat \Gamma}(\hat u)$, as an orbit of $u$
is $\ell_{\hat \Gamma}(\hat u)$-dense on $\alpha$.  Let $g' = u^n g$.  Then $g'$
serves to conjugate $u$ to $v$.  Moreover,
$d_{\hat \Gamma}(e, g') \le d_{\hat \Gamma}(e, a) + d_{\hat \Gamma}(a, g'b) + d_{\hat \Gamma}(g'b, g')$.
Thus $\ell_{\hat \Gamma}(\hat g') = d_{\hat \Gamma}(e, g') \le K \ell_{\hat \Gamma}(\hat u) + 
\ell_{\hat \Gamma}(\hat u) + K \ell_{\hat \Gamma}(\hat v)$.
The lemma follows by observing that $\ell_{\Gamma}(u) \ge \ell_{\hat \Gamma}(\hat u)$ and 
$\ell_{\Gamma}(v) \ge \ell_{\hat \Gamma}(\hat v)$.
\hfill //
\endpf

\proclaim{\bf\underbar{Lemma 2.2.9}}
Let $u$ be a nontorsion parabolic elements of $G$, lying in $\gamma H \gamma^{-1}$.  
Then there exists a universal polynomial $K_p$, of degree one, such that
$\ell_{\hat \Gamma}(\hat \gamma) \le K_p ( \ell_{\Gamma}( u ) )$.
\endproclaim
\prf
For a coset $\gamma H$, let 
$\ell^*_{\Gamma}(\gamma H) = \min\{ \ell_{\Gamma}(xh) \text{ for } h \in H\}$.
Let $x$ be a point of $\gamma H$ such that $\ell_{\Gamma}(x) = \ell^*_{\Gamma}(\gamma H)$.
By Lemma 4.24 of [DS] there exists a constant $C$, independant of $u$, $\gamma$ and $H$,  
such that either $d_{\Gamma}(x,ux) \le C$ or 
$\ell_{\Gamma}(u) \ge \ell^*_{\Gamma}(\gamma H) + \frac{1}{2} d_{\Gamma}(x, ux) - C$.

Thus if $d_{\Gamma}(x,ux) > C$ then 
$\ell^*_{\Gamma}(\gamma H) \leq \ell_{\Gamma}(u) - \frac{1}{2} d_{\Gamma}(x, ux) + C$
As $\ell_{\hat \Gamma}(\hat \gamma) \le \ell^*_{\hat \Gamma}(\gamma H) + 1 \le 
	\ell^*_{\Gamma}(\gamma H) + 1$, it follows that
$\ell_{\hat \Gamma}(\hat \gamma) \le \ell_{\Gamma}(u) + C + 1$.

Otherwise $\ell_{\Gamma}(x^{-1}ux) = d_{\Gamma}(x,ux) \le C$.  Consider the sequence
of elements $x^{-1} u^n x \in G$.  As $u$ is nontorsion, this is an infinite family
of elements.  As $\ell_{\Gamma}$ is a proper length function, there can be at most
finitely many elements of $G$ with length not exceeding $C$.  Let $M$ be the
cardinality of the ball of radius $C$ centered at the identity.  Then there is
an $N \le M+1$ such that $\ell_{\Gamma}(x^{-1} u^N x) > C$.
As $u^N$ is a nontorsion parabolic element lying in $\gamma H \gamma^{-1}$, one has
$\ell_{\Gamma}(u^N) \ge \ell^*_{\Gamma}(\gamma H) + \frac{1}{2} d_{\Gamma}(x, u^N x) - C$
so that $\ell_{\hat \Gamma}(\hat \gamma) \le \ell_{\Gamma}(u^N) + C + 1$.
As $\ell_{\Gamma}(u^N) \le N \ell_{\Gamma}(u) \le (M+1)\ell_{\Gamma}(u)$, we have
$\ell_{\hat \Gamma}(\hat \gamma) \le (M+1)\ell_{\Gamma}(u) + C + 1$.

\hfill //
\endpf

\proclaim{\bf\underbar{Theorem 2.2.10}} Let $G$ be a finitely
generated group with generating set $X$, such that $G$ is
relatively hyperbolic with respect to $\{ H_\lambda \}_{\lambda
\in \Lambda}$, each member of which has a $\Cal P$-solvable conjugacy-bound.
Then $G$ has a $\Cal P$-solvable conjugacy-bound.  
\endproclaim

\prf
Let $\delta$ denote the hyperbolicity constant of the coned-off Cayley graph $\hat \Gamma$.
As each $H_\lambda$ has a $\Cal P$-solvable conjugacy-bound, there are 
polynomials $P_\lambda$ such that if $u$, $v \in H_\lambda$ are conjugate
in $H_\lambda$, then there is a $g \in H_\lambda$ with $u^g = v$ and
$\ell_{H_\lambda}( g ) \le P_{\lambda} (\ell_{H_\lambda}( u ) + \ell_{H_\lambda}( v ))$.
As $(H_\lambda, \ell_{H_\lambda}) \to (G, \ell_G)$ is a quasi-isometric
embedding, there exists a polynomial $Q_\lambda$ such that
$\ell_\Gamma (g) \le Q_{\lambda}( \ell_\Gamma (u) + \ell_\Gamma(v) )$.
If $u$ and $v$ are two elements of $H_\lambda$, with $g \in G$ with
$u^g = v$, then $v \in H_\lambda^g \cap H_\lambda$.  If 
$g \in G \setminus H_\lambda$, then Proposition 2.36 of [Os] gives $v=1$.
Thus if $u$ and $v$ are nontrivial elements of $H_\lambda$ which have a
conjugator in $G$, then they have a conjugator in $H_\lambda$.  The $\Gamma(G,X)$-length
of the minimal conjugator is thus bounded by the polynomial $Q_\lambda$.
Moreover, as $\Lambda$ must be finite by Proposition 2.2.1 above, 
$Q = \sum_{\lambda \in \Lambda} Q_\lambda$ is a polynomial which works to bound 
the conjugator length for every $H_\lambda$.

Let $u$ and $v$ be conjugate hyperbolic elements of $G$, and let
$L = \ell_\Gamma(u) + \ell_\Gamma(v)$
According to [Bu] Lemma 5.11 we have that the $\Gamma$-distance which
$g$ travels in a coset it penetrates is bounded by $2( L ) + 10 c( 8L )$,
where $c(k)$ is the function given by the BCP property.  As shown above
this is bounded by a polynomial in $L$.
By Lemma 2.2.8 above, $\ell_{\hat \Gamma}( \hat g ) \le K_h( L )$.
As each step of this shortened path corresponds to a $\Gamma$-distance
bounded by a polynomial in $L$, it follows that there is an element $g' \in G$
conjugating $u$ to $v$ with
$\ell_\Gamma(g') \le \left(2( \ell_\Gamma(u) + \ell_\Gamma(v) ) + 
	10 c( 8 \ell_\Gamma(u) + 8 \ell_\Gamma(v) )\right)K_h( \ell_\Gamma(u) + \ell_\Gamma(v) )$.
As $c(k)$ can be bounded by a polynomial in $k$, we see that this expression
is bounded by some universal polynomial $Q_h( \ell_\Gamma(u) + \ell_\Gamma(v) )$.

Let $u$ be a nontorsion parabolic element of $G \setminus H_\lambda$,
and $v$ an element of $H_\lambda$ with $u$ conjugate to $v$.  Let $g$ 
be a conjugator of minimal relative length with $g^{-1}ug = v$.  Then 
$u \in g H_\lambda g^{-1}$.  By Lemma 2.2.9 we have that
$\ell_{\hat \Gamma}(\hat g) \le K_p ( \ell_{\Gamma}( u ) )$.

Let $g'$ be an element with minimal relative length which conjugates $u$ into $H$, and let
$h = {g'}^{-1}ug'$.  By Theorem 5.13 of [Bu] we have the $\Gamma$-length is bounded by
$c( 7\ell_\Gamma(u) + 7\ell_\Gamma(v) )$.  In this case, from our above argument, we have that
the element of minimal length conjugating $h$ to $v$ has length bounded by a polynomial in
$\ell_\Gamma(u) + \ell_\Gamma(v)$.  Then to bound the length of $g$, we may assume that
$g$ has minimal relative length amongst all elements conjugating $u$ into $H$.

Consider the geodesic quadrilateral in $\hat \Gamma$ with sides
$[e,u]$, $[u,ug]$, $[ug, g]$ and $[e,g]$, where the paths $[e,g]$ and $[ug,g]$
correspond to the shortened path $\hat g$.  Denote by $\hat p$ the path $[u,ug]$
and by $\hat q$ the path $[e,g]$.  Let $\hat p$ travel through the coset $fH$ along
the word $k$.  By [Bu] Corollary 5.4, if $\hat q$ 
do not penetrate $fH$, then $\Gamma$-distance of $k$ is at most 
$2\left( \ell_\Gamma(u) + \ell_\Gamma(v) \right) + 	
	c\left( 2 \ell_\Gamma(u) + 2 \ell_\Gamma(v) + 1 \right) + 2c(2)$.  
Otherwise $\hat q$ also travels through $fH$.  Assume that $\hat q$ travels
through $fH$ along $k'$.

If $k = k'$ so that $p$ and $q$ travel synchronously through $fH$, let $\gamma_1$ be
an $H$-geodesic connecting the first point of $p$ in $fH$ to the first point of $q$
in $fH$.  Let $\gamma_2$ be an $H$-geodesic connecting the last point of $p$ in $fH$ 
to the last point of $q$ in $fH$.  Then $u$ and $v$ are conjugate to $\gamma_1$ and
$\gamma_2$, each of which are in $H$.  By the argument in the proof of Theorem 5.12
of [Bu], the relative lengths of $\gamma_1$ and $\gamma_2$ are each bounded by
$$
\left( 7\ell_\Gamma(u) + 7\ell_\Gamma(v) + 2 \delta + 
	2 \tilde N(2\ell_\Gamma(u)+2\ell_\Gamma(v) + 1)\right)
$$
Moreover, in any coset they penetrate	they can travel a $\Gamma$-distance at most
$$
	\left(  3 c( 56\ell_\Gamma(u) + 56\ell_\Gamma(v) + 18 \delta + 
	16 \tilde N(2\ell_\Gamma(u)+2\ell_\Gamma(v) + 1) ) + 
	4\ell_\Gamma(u) + 4\ell_\Gamma(v) + 10c(8\ell_\Gamma(u) + 8\ell_\Gamma(v))\right)
$$
Thus their $\Gamma$-lengths are bounded by the product of these two quantities, a polynomial
in $\ell_\Gamma(u) + \ell_\Gamma(v)$.
The subpath of $g$ conjugating $\gamma_1$ to $\gamma_2$ can then be thought of as an element
of $H$ conjugating $\gamma_1$ to $\gamma_2$.  It's $\Gamma$-length can then be
bounded by a polynomial in $\ell_\Gamma( \gamma_1) + \ell_\Gamma( \gamma_2 )$.
Thus the $\Gamma$-distance traveled by $p$ and $q$ through $fH$ is bounded by
a universal polynomial in $\ell_\Gamma(u) + \ell_\Gamma(v)$.

If this is not the case, then $k$ and $k'$ are not synchronously in $fH$.
Assume that there is a coset $f'H$ penetrated by both $\hat p$ and $\hat q$, such
that $\hat p$ travels through $f'H$ along $j$, $\hat q$ travels through $f'H$ along $j'$
with $\hat p$ traveling along $k$ before $j$ and $\hat q$ traveling
along $j'$ before $k'$.  

\centerline{\epsfbox{skew.eps}}


Let $a$ be the first point of $q$ in $f'H$, and let $b$ be the first point of $p$ in $f'H$, and
consider the paths $w_1$ and $w_2$ where $w_1$ is the indicated paths from $e$ to $u$ then to $a$,
and $b$ is the path from $e$ to $b$ then to $a$, such that the path from $b$ to $a$ is an
$H$-geodesic.  As in the `skew'-coset case of Lemma 5.5 of [Bu], $\hat w_1$ and $\hat w_2$ are two
$\ell_\Gamma(u) + \ell_\Gamma(v) + 2 \ell_{\hat \Gamma}(\hat g)$ - quasigeodesics, and
the $\Gamma$-length of $k$ is bounded by 
$c( \ell_\Gamma(u) + \ell_\Gamma(v) + 2 \ell_{\hat \Gamma}(\hat g) )$.  
As $\ell_{\hat \Gamma}(\hat g)$ is bounded by a polynomial in $\ell_\Gamma(u) + \ell_\Gamma(v)$,
we have that the $\Gamma$-length of $k$ itself is bounded by a polynomial in 
$\ell_\Gamma(u) + \ell_\Gamma(v)$.

If this is not the case, then we perform a series of surgeries on the paths $p$ and $q$ as follows.
Let $I$ denote the set of $i$ such that $p(i)$ is a point where $p$ enters or exits a coset $f'H$ which 
both $p$ and $q$ penetrate, as well as $i$ such that $q(i)$ is a point where $q$ enters or exits a coset 
$f'H$ which both $p$ and $q$ penetrate.
Connect $p(i)$ to $q(i)$ by a $\Gamma$-geodesic, denoted by $\gamma_i$, for each $i \in I$.
The $\Gamma$-length of each $\gamma_i$ is bounded in the same way as $\gamma_1$ and $\gamma_2$ above.
This splits the paths $p$ and $q$ into a series of subpaths which conjugate one $\gamma_i$
to the next.  We will denote by $\gamma_0$ and $\gamma_\infty$ the paths $u$ and $v$ respectivly,
and consider $0$ and $\infty$ as elements of $I$.

\centerline{\epsfbox{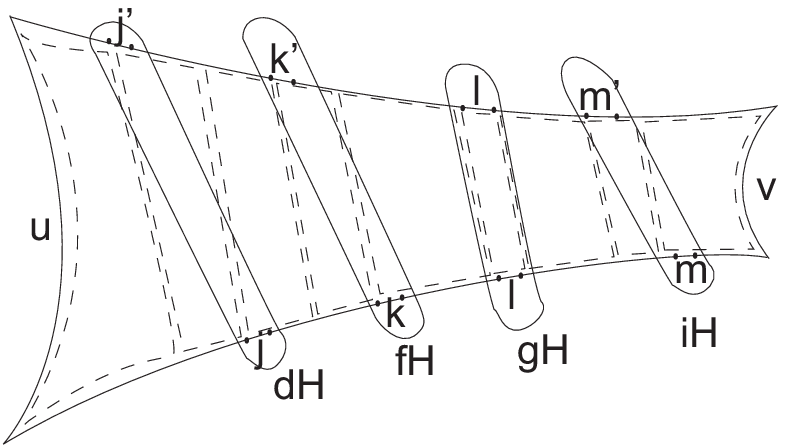}}

A priori, for each successive $\gamma_i$, $\gamma_j$ conjugacy, we have one of the following:

(1) $\gamma_i$ and $\gamma_j$ are in the same coset.

(2) $\gamma_i$ and $\gamma_j$ do not lie in a single common coset, and the subpaths $p'$ and $q'$ connecting 
    $\gamma_i$ to $\gamma_j$ do not penetrate a common coset.

(3) $\gamma_i$ and $\gamma_j$ do not lie in a single common coset, and the subpaths $p'$ and $q'$ connecting 
    $\gamma_i$ to $\gamma_j$ do penetrate a common coset.
    
In each of these cases, the subword of $g$ which conjugates $\gamma_i$ to $\gamma_j$ can have its
length bounded by that of the minimal length conjugator of $\gamma_i$ and $\gamma_j$.
In case (1) the conjugator is a word in $H$ which we have already seen has
length bounded by a polynomial in the $\Gamma$-lengths of $\gamma_i$ and $\gamma_j$.  
Thus this case  is solvable.
In case (2) the argument above where the paths $p$ and $q$ do not penetrate a common coset 
bounds the length of the minimal conjugator by a polynomial in the $\Gamma$-lengths of 
$\gamma_i$ and $\gamma_j$.  Thus this case is also solvable.
In case (3), the subpaths $p'$ and $q'$ of $p$ and $q$ connecting $\gamma_i$ to $\gamma_j$
both travel a positive distance through a common coset.  Denote this coset by $fH$.  By the
definition of $I$, $p(i)$ is the point where $p$ enters $fH$ and $p(j)$ is the point where
$p$ exits $fH$, and similarly for $q(i)$ and $q(j)$.  Then $\gamma_i$ and $\gamma_j$ are
elements of $H$, and the above argument bounds the length of their conjugator.
 
As the number of such pairs is bounded by the $\hat \Gamma$-length of $\hat g$, which is itself
bounded by a polynomial in $\ell_\Gamma(u) + \ell_\Gamma(v)$, and each $\gamma_i$ has length
bounded by a polynomial in $\ell_\Gamma(u) + \ell_\Gamma(v)$.
Putting this together with the above arguments, we see that, through any coset which they 
penetrate, the $\Gamma$-distance which $p$ and $q$ can travel is bounded by a polynomial
in $\ell_\Gamma(u) + \ell_\Gamma(v)$.
Thus there is a universal polynomial $U$ such that the element of minimal length conjugating
$u$ to $v$ has length bounded by $U( \ell_\Gamma(u) + \ell_\Gamma(v) )$.

Let $u$ and $v$ be two conjugate parabolic elements of $G \setminus H_\lambda$,
which can each be conjugated into $H_\lambda$.
Assume that $u^g = v$, for $g$ a conjugator of minimal length.  
There is $h \in H_\lambda$ such that $u$ is conjugate to $h$, and $v$ is
conjugate to $h$, and by Lemma 5.13 of [Bu] we have 
$\ell_\Gamma( h ) \le c(7 \ell_\Gamma(u) + 7 \ell_\Gamma(v) )$.
By the argument above, there is an element conjugating
$u$ to $h$ whose length is bounded by a polynomial in $\ell_\Gamma( u ) + \ell_\Gamma( h )$
and an element conjugating $v$ to $h$ whose length is bounded by a polynomial
in $\ell_\Gamma( v ) + \ell_\Gamma( h )$.  The product of these two conjugators gives an element
of $G$ which conjugates $u$ to $v$, and whose length is bounded by a universal
polynomial in $\ell_\Gamma(u) + \ell_\Gamma(v)$.
\hfill //
\endpf

\exsp Recall that if $L_1$ and $L_2$ are two length
functions on a group $G$, then  $L_1$ dominates
$L_2$ if there are constants $A$ and $B > 0$ such that for all
$g \in G$ $L_2(g) \le A(L_1(g) + 1)^B$.  $L_1$ and
$L_2$ are equivalent if $L_1$ dominates $L_2$ and
$L_2$ dominates $L_1$.  Viewed as weight functions on $G$, we see 
that this classical notion of equivalence corresponds to 
$\Cal P$-equivalence as defined in section 1.2 above. Thus 
equivalent length functions on a discrete group yield the same 
algebra of polynomially bounded functions.

\exsp \proclaim{\bf\underbar{Theorem 2.2.11}} Let $G$ be a finitely
generated group, such that $G$ is relatively hyperbolic with
respect to $H_1, H_2, ... , H_n$, where each $H_i$ is in the class
$\cC_{\Cal P}(\Bbb C)$ and has a $\Cal P$-solvable conjugacy-bound. 
Then $G$ is in the class $\cC_{\Cal P}(\Bbb C)$.
In particular, the $\ell^1-$SrBC holds for $G$. 
\endproclaim

\prf By Theorem 2.2.10 and the methods of section 1.4 above,
$$
\gather
HH_t^*(H^{1,\infty}_L(G))_x \cong HH_t^*(H^{1,\infty}_L(G_h))\\
HC_t^*(H^{1,\infty}_L(G))_x \cong HC^*_t(H^{1,\infty}_L(G_h))
\endgather
$$

It remains to show that $G_h$ belongs to the class
$\cC_{\Cal P}(\Bbb C)$, for each non-torsion element $h$. 

\exsp Let $h \in G$ be a nontorsion element.

\exsp If $h$ is hyperbolic, then by Prop. 2.2.2, the embedding
$G_h \to G$ is a quasi-isometric embedding.  This is precisely the
statement that $L_{G}|_{G_h}$ is equivalent to $L_{G_h}$, where $L_G$ 
denotes the standard word-length function on $G$ (determined by the 
generating set). We also have $G_h$ is virtually cyclic. Thus, 
${G_h} \in \cC_{\Cal P}(\Bbb C)$.

\exsp Suppose that $h$ is parabolic. Then by Prop. 2.2.3, 
$G_h = x H_{i, x^{-1}hx} x^{-1}$, for some $x \in G$. 
Thus, $L_{G_h} = L_{x (H_{i,
x^{-1}hx} )x^{-1}}$.  As all of the $H_i$ are quasi-convex in $G$,
the restriction $L_G |_{x H_i x^{-1}}$ is equivalent to
$L_{x H_i x^{-1}}$.  Then $L_G | _{x H_{i, x^{-1}hx}
x^{-1}}$ is equivalent to $L_{x H_i x^{-1}}|_{x H_{i, x^{-1}hx}
x^{-1}}$. That is, $L_G|_{G_h}$ is equivalent to $L_{x H_i
x^{-1}}|_{x( H_{i, x^{-1}hx} )x^{-1}}$, which is by hypothesis
equivalent to $L_{x H_{i,h} x^{-1}}$.  Conjugating where
appropriate we obtain $L_G|_{G_h}$ is equivalent to
$L_{G_h}$.  Since $G_h \cong H_{i,
x^{-1}hx}$, while $H_i$ are in the class $\cC_{\Cal P}(\Bbb C)$ the theorem is proved.\hfill //
\endpf
\vskip.2in

This verifies Theorem A in the case $\Cal B = \Cal P$. However, in the more general case when $\Cal P\prec \Cal B$, the argument is the same. In fact, the $\Cal B$-analogue of Theorem 2.2.10 clearly holds whenever $\Cal P\prec \Cal B$, at which point the argument for proving Theorem 2.2.11 carries over with $\Cal B$ replacing $\Cal P$ throughout.
\vskip.3in

\subhead 2.3. Combable groups
\endsubhead
\vskip.1in

We begin with the definition of a \underbar{combing} for a
discrete group given in [Me1]( this notion is originally due to
Thurston [CEHLPT]).

\proclaim{\bf\underbar{Definition 2.3.1}} If $(X,*)$ is a
basepointed discrete metric space with distance function $d$, a
\underbar{combing} on $X$ consists of a sequence of maps $f_n:X\to
X$ satisfying \vskip.05in

(i) $f_0 = *$ and for all $x\in X$ there is an $N$ such that
$f_n(x) = x\ \forall n\ge N$; \vskip.05in

(ii) there exists a constant $C > 0$ such that
$d(f_n(x),f_n(y))\le Cd(x,y) + 1$; \vskip.05in

(iii) there exists a constant $S\in\Bbb N$ such that
$d(f_n(x),f_{n+1}(x)) \le S$ for all $x\in X$ and $n\in\Bbb N$.
\endproclaim

For a given $x$, let $J(x)$ be the number of integers $n$ for
which $f_n(x)\ne f_{n+1}(x)$. The combing is \underbar{polynomial}
resp. \underbar{exponential} if  $J(x)\le g(l(x))$ for all $x\in
X$, where $g(t) = C'(1+t)^m$ resp. $C'e^{kt}$ for some
$C',k\in\Bbb R_+$ and $m\in\Bbb N$ (here $l(x) := d(*,x)$). This
applies in particular to discrete groups with word-length $(G,L)$,
viewed as a basepointed metric space with metric $d(g_1,g_2) :=
L(g_1^{-1}g_2)$.

That $HH^*_t(H_{\Cal B,L}(G))$ and $HC^*_t(H_{\Cal B,L}(G))$
should be computable in terms of the Hochschild and cyclic
cohomology of $\Bbb C[G]$ for combable groups is strongly
suggested by the following two results.

\proclaim{\bf\underbar{Theorem 2.3.2 - [Me1], [Me2], [Og2]}} Let
$(G,L)$ be a synchronously combable group with polynomial combing
function, $L$ the standard word-length function on $G$. Then there
is an isomorphism
$$
H^*_{\Cal P,G}(EG.;A)\overset\cong\to\longrightarrow H^*_G(EG.;A)
$$
for any p-bounded semi-normed (p.s.) $G$-module $A$.
\endproclaim

Actually the result of [Me1] only deals with the case $A = $ the
trivial $G$-module $\Bbb C$, but also proves the same result when
$G$ has an exponential combing function and $\Cal P$ is replaced
by $\Cal E$; for more general coefficient modules one needs to use
either [Me2] or [Og2] (the definition of a p.s. $G$-module is as in
[Og2]).

\proclaim{\bf\underbar{Theorem 2.3.3 - [Me2]}} Let $(G,L)$ be as
above, and let $P_*$ be a bimodule resolution of $\Bbb C[G]$ by
free $\Bbb C[G]$-bimodules. Then $H_{\Cal B,L}(G)\underset{\Bbb
C[G]}\to\otimes P_*\underset{\Bbb C[G]}\to\otimes H_{\Cal B,L}(G)$
is a (bornological) bimodule resolution of $H_{\Cal B,L}(G)$ by
free (bornological) $H_{\Cal B,L}(G)$-bimodules ($\Cal B = \Cal P$
or $\Cal E$).
\endproclaim

[Note: What is done in [Me2] is much more general; Theorem 2.3.2
only states a special case relevant to our setting].

\proclaim{\bf\underbar{Proposition 2.3.4}} Assume that $(G,L)$ is
a group equipped with a polynomial combing in the above sense.
Then for each conjugacy class $x$ there are isomorphisms
$$
HH_t^*(H_{\Cal P,L}(G))_{x} = HH^*_{\Cal P}(\Bbb C[G])_{x} =
H^*_{\Cal P,G}(EG.;Hom_{\Cal P}(S_x,\Bbb C))\cong
H^*_{G}(EG.;Hom_{\Cal P}(S_x,\Bbb C))
$$
where (i) $S_x$ is a weighted $G$-set with weight function induced
by $L$ and the inclusion into $G$, (ii) $Hom_{\Cal P}(S_x,\Bbb C)$
denotes the subspace of $Hom(S_x,\Bbb C)$ consisting of those maps
which are polynomially bounded with respect to the weight function
on $S_x$, and (iii) $H^*_{\Cal P,G}(_-)$ denotes the cohomology of
the cocomplex of $G$-equivariant, polynomially bounded cochains.
\endproclaim

\prf The first equality was verified above. Corollary 1.4.2 gives
the identification of $HH^*_{\Cal P}(\Bbb C[G])_{x}$ with
$H^*_{\Cal P}( S_x\underset{G}\to\times EG.) = H^*_{\Cal P,G}(
S_x\times EG.)$. These last groups are computed as the cohomology
of $\{C^n_{\Cal P,G}(S_x\times EG.,\Bbb C)\}_{n\ge 0} =
\{Hom_{\Cal P,G}(S_x\times EG_n,\Bbb C)\}_{n\ge 0}$. For each $n$
taking partial adjoints yields an isomorphism $Hom_{\Cal
P,G}(S_x\times EG_n,\Bbb C)\cong Hom_{\Cal P,G}(EG_n,Hom_{\Cal
P}(S_x,\Bbb C))$ - this is just the usual isomorphism on cochains,
restricted to polynomially bounded cochains (cf. [Og2]);
consequently, it commutes with the coboundary maps and extends to
an isomorphism of cocomplexes
$$
\{Hom_{\Cal P,G}(S_x\times EG_n,\Bbb C)\}_{n\ge 0}\cong C^*_{\Cal
P,G}(EG.,Hom_{\Cal P}(S_x,\Bbb C)):= \{Hom_{\Cal
P,G}(EG_n,Hom_{\Cal P}(S_x,\Bbb C))\}_{n\ge 0}
$$
which produces the second isomorphism in cohomology. The third
follows from Theorem 2.3.2.\hfill //
\endpf

Note that $Hom_{\Cal P}(S_x,\Bbb C)$ is dense in $Hom(S_x,\Bbb
C)$, and it is reasonable to ask (or even conjecture)

\proclaim{\bf\underbar{Question/Conjecture 2.3.5}} For $(G,L)$ as
above, does the inclusion of coefficients $Hom_{\Cal P}(S_x,\Bbb
C)\hookrightarrow Hom(S_x,\Bbb C)$ induce an isomorphism
$H^*_{G}(EG.;Hom_{\Cal P}(S_x,\Bbb
C))\overset\cong\to\longrightarrow H^*_{G}(EG.;Hom(S_x,\Bbb C))$?
\endproclaim

If this is the case, then for each $x$ there would be an
isomorphism $HH_t^*(H_{\Cal P,L}(G))_{x}\cong HH^*(\Bbb
C[G])_{x}$, and consequently isomorphisms $HC_t^*(H_{\Cal
P,L}(G))_{x}\cong HC^*(\Bbb C[G])_{x}$. However, we have not been
able to verify these isomorphisms except in certain special cases
in which the general theory of section 1.4 can be applied. The
following is an immediate consequence of
Proposition 2.3.4.

\proclaim{\bf\underbar{Corollary 2.3.6}} If $\Bbb Q$ admits a
projective resolution of length $N$ over $\Bbb Q[G]$, then
$HH_t^*(H_{\Cal P,L}(G))_{x} = HH^*(\Bbb C[G])_{x} = 0$ for all 
$* > N$.
\endproclaim

An important subclass of combable groups are the
\underbar{semi-hyperbolic groups} defined by Alonso and Bridson
in their paper [AB]. We have

\proclaim{\bf\underbar{Theorem 2.3.7}} Let $(G,L)$ be
semi-hyperbolic in the sense of [AB]. Then for each conjugacy
class $x=<h>$, there are isomorphisms $HH_t^*(H_{\Cal
E,L}(G))_{x}\cong H^*(BG_h;\Bbb C) = HH^*(\Bbb C[G])_{x}$,
$HC_t^*(H_{\Cal E,L}(G))_{x}\cong HC^*(\Bbb C[G])_{x}$.
\endproclaim

\prf By [AB], we know that quasi-convex subgroups of
semi-hyperbolic are again semi-hyperbolic. This applies in
particular to the centralizer subgroups $G_x$. Consequently, each
$G_x$ admits a synchronous (linear) combing. By [Me1], $G_x$ is
$\Cal E$-isocohomogical. Again, by [AB], we have that $G$ has an
 $\Cal E$-solvable conjugacy-bound. Thus the hypothesis of
Theorem 1.4.7 is satisfied for any $\Cal B$ containing $\Cal E$,
and the result follows.\hfill //
\endpf

\proclaim{\bf\underbar{Corollary 2.3.8}} If $(G,L)$ is
semi-hyperbolic and $G$ satisfies the nilpotency condition, then
it satisfies the $\Cal E$-nilpotency condition.
\endproclaim

In the absence of an affirmative answer to (2.3.5) above, it would
be quite desirable to strengthen Corollary 2.3.8 by replacing
\lq\lq $\Cal E$-nilpotency\rq\rq with \lq\lq $\Cal
P$-nilpotency\rq\rq. Certainly the groups $G_x$ are $\Cal
P$-isocohomological (their combings are even linear). The
obstruction lies with the time needed to solve the conjugacy
problem. Currently it is not known if the exponential bound on the
solution time given in [AB] is actually a sharp bound in the
general case, or whether it can be improved to polynomial time.
However, for certain semi-hyperbolic groups, this improved bound
is known for certain conjugacy classes. In particular, if $G$ is a mapping class group, then the
conjugacy-bound is polynomial [MM] when the conjugacy class is pseudo-Anosov.
In general the ${\Cal P}$-solvable conjugacy-bound problem is still solvable but we will address this issue in a future publication.

\bigskip
Finally we remark that the Bass Conjecture for semi-hyperbolic groups was proved by Eckman [Ec2]

\vskip.3in

\subhead 2.4. Hochschild and cyclic homology of $\ell^1(G)$
\endsubhead
\vskip.1in

As we have observed, for $\Cal B = \Cal B_{min}$, the
corresponding rapid decay algebra $H_{\Cal B_{min},L}(G)$ is equal to
$\ell^1(G)$, the full $\ell^1$-algebra of $G$. More generally, one has an
identification of Banach spaces $H_{\Cal B_{min},w}(X)=\ell^1(X)$
for any weighted space $(X,w)$. Since the class $\Cal B_{min}$
consists of constant functions, the space $H_{\Cal B_{min},w}(X)$
is independent of $w$. For this reason, we will drop $w$ and
simply write $\ell^1(X)$ for the corresponding space. As in the
proof of Lemma 1.3.3, there is, for any two sets $X, X'$, a
natural isomorphism of Banach spaces
$$
\ell^1(X)\tensor \ell^1(X')\cong \ell^1(X\times X')
\tag2.4.1
$$
and by Lemma 1.2.4 a natural identification
$$
\ell^1(X)^*\cong Hom_b(X,\Bbb C)
\tag2.4.2
$$
between the continuous dual of $\ell^1(X)$ and the space of
bounded, complex-valued maps on $X$. Consequently, as in Lemma
1.3.6, we have an isomorphism
$$
C^*_t(\ell^1(G))\cong C^*_b(N.^{cy}(G))
$$
between the topological Hochschild cocomplex for $\ell^1(G)$ and
the cocomplex of bounded cochains on $N.^{cy}(G)$. Again, this
isomorphism preserves the respective decompositions by conjugacy
classes, and commutes with the maps in the Connes-Gysin seqence.
Thus for any conjugacy class $x$ we have an identification
$$
HH^*_t(\ell^1(G))_x\cong H^*_b(N.^{cy}(G)_x)
$$
Unlike the case $\Cal B\ne \Cal B_{min}$, we see that the
isomorphism $G_h\backslash G\cong S_x$ of right $G$-sets
(Proposition 1.4.3, Lemma 1.4.4) producing the isomorphism of
simplicial sets $(\pi_x).: G_h\backslash G \underset{G}\to\times
EG.\overset\cong\to{\longrightarrow} S_x \underset{G}\to\times
EG.$, will always yield an isomorphism on bounded cochains:
$C^*_b(S_x \underset{G}\to\times
EG.)\overset\cong\to\longrightarrow C^*_b(G_h\backslash G
\underset{G}\to\times EG.)$. Given that Proposition 1.4.5 holds
without restriction, we conclude

\proclaim{\bf\underbar{Corollary 2.4.3}} For all discrete groups
$G$ and conjugacy classes $x=<h>$, there is a natural isomorphism
$HH^*_t(\ell^1(G))_x\cong H^*_b(BG_h)$.
\endproclaim

We will write $S^*_{\Bbb C}$ for the $S$-map $S^*:HC^*(\Bbb C)\to
HC^{*+2}(\Bbb C)$. The main result of this section is

\proclaim{\bf\underbar{Theorem 2.4.4}} For all discrete groups $G$
and conjugacy classes $x=<h>$ there is an isomorphism
$$
HC^*_t(\ell^1(G))_x\cong H^*_b(BG_h)\otimes HC^*(\Bbb C)
$$
where the map $S^*_x$ on the left identifies with $Id\otimes S^*_{\Bbb C}$ on the right.
\endproclaim

\prf The proof is an adaption of the method due to Nistor [Ni1](see also [Ji2]).
For $h\in G$, let $\wt{L}(G,h). = (EG.)$. Writing $x$ for $<h>$, define a projection map
$P_h:\wt{L}(G,h).\surj N.^{cy}(G)_x$ by $P_h(g_0,\dots,g_n) =
(g_n^{-1}hg_0\otimes g_0^{-1}g_1\otimes\dots\otimes
g_{n-1}^{-1}g_n)$. It is easy to see that $P_h$ induces an
isomorphism of simplicial sets $\ov{P_h}:*\underset G_h\to\times
\wt{L}(G,h).\overset\cong\to\longrightarrow N.^{cy}(G)_x$, and
hence an isomorphism of bounded cochain cocomplexes
$$
C_b^*(N.^{cy}(G)_x)\overset\cong\to\longrightarrow
C_b^*(*\underset G_h\to\times \wt{L}(G,h).)\cong
C_{b,G_h}^*(\wt{L}(G,h).)
$$
where $C_{b,G_h}^*(\wt{L}(G,h).)$ denotes the subcocomplex of
$C_{b}^*(\wt{L}(G,h).)$ consisting of $G_h$-equivariant
cochains (with coefficients in the trivial $G_h$-module $\Bbb C$).
Now consider the map $T\hskip-.03in.$ defined on $\wt{L}(G,h).$ by
$T_{n+1}(g_0,\dots,g_n) = (hg_n,g_0,\dots,g_{n-1})$. If $H$ is any
subgroup of $G_h$ containing $h$, then we see that $T\hskip-.03in.$ descends
to a map $\ov T\hskip-.03in.$ on $*\underset H\to\times \wt{L}(G,h).$,
giving $*\underset H\to\times \wt{L}(G,h).$ the structure of a
cyclic simplicial set, and making $\ov{P_h}$ an isomorphism of cyclic
simplicial sets. Denote the singular resp. cyclic complex of a
cyclic simplicial set $X.^{c}$ (with coeff.s in $\Bbb C$) by
$C_*(X.^{c})$ resp. $CC_*(X.^{c})$, and the corresponding dual
cocomplexes by $C^*(X.^{c})$ resp. $CC^*(X.^{c})$. Clearly an
isomorphism $\phi:X.^c\overset\cong\to\to Y.^c$ of cyclic
simplicial sets induces an isomorphism on both singular and cyclic
(co)complexes. It also induces an isomorphism of
\underbar{bounded} cyclic cocomplexes
$CC_b^*(\phi):CC^*_b(X.^c)\overset\cong\to\to CC^*_b(Y.^c)$ (as we
are in char. 0, $CC^*_b(Z.^c)$ is quasi-isomorphic to the
subcocomplex of $C^*_b(Z.^c)$ consisting of cochains invariant
under the cyclic operator). Applying these observations to our
situation, we get for each conjugacy class $x = <h>$ an
isomorphism
$$
CC_t^*(\ell^1(G))_x\cong CC^*_{b,G_h}(\wt{L}(G,h).)
$$
Let $C_h\subset G_h$ denote the cyclic subgroup generated by $h$,
and $N_h$ the quotient of $G_h$ by $C_h$. Denote the cohomology
groups of  $CC^*_{b,C_h}(\wt{L}(G,h).)$ by
$HC^*_{b,C_h}(\wt{L}(G,h).)$. We claim these groups are
Hausdorff, or alternatively, that the reduced and unreduced
cohomology groups agree (see Note (1.3.2)). In fact, as $C_h$ is amenable for all
$h$ (it is abelian), the augmentation map $\wt{L}(G,h).\surj *$ induces
isomorphisms on bounded Hochschild cohomology, and hence also bounded cyclic cohomology.
The Hausdorff condition guarantees the existence of a
first-quadrant spectral sequence strongly converging to
$HC^*_{b,G_h}(\wt{L}(G,h).):= H^*(CC^*_{b,G_h}(\wt{L}(G,h).))$
with $E_2^{**}$-term
$$
E_2^{**} = H_b^*(N_h;HC^*_{b,C_h}(\wt{L}(G,h).))
$$
(compare [No], [Ji1, (3.7)]). As we have just observed, the
augmentation map on $\wt{L}(G,h).$ yields an isomorphism
$HC^*_{b,C_h}(\wt{L}(G,h).)\overset\cong\to\longrightarrow
HC^*(\Bbb C)$, and hence an isomorphism
$$
E_2^{**} = H_b^*(N_h;HC^*_{b,C_h}(\wt{L}(G,h).))
\overset\cong\to\longrightarrow H_b^*(N_h;HC^*(\Bbb C))
$$
of $E_2^{**}$-terms. The spectral sequence on the right converges
to $H^*_b(N_h)\otimes HC^*(\Bbb C)$; however by the Serre spectral
sequence in bounded cohomology associated to the exact sequence
$C_h\inj G_h\overset{proj_h}\to\surj N_h$ (loc. sit.), we have
$H^*_b(N_h)\overset{proj_h}\to{\underset\cong\to\longrightarrow}
H^*_b(G_h)$, completing the proof.\hfill
//
\endpf

Thus, although computable, the non-elliptic summands of the
topological cyclic cohomology groups of $\ell^1(G)$ will never
have a nilpotent $S$-operator, and so are not useful in verifying
the $\ell^1$-Bass conjecture by the methods used in this paper.
Note that if $G$ itself is amenable, then each $G_x$ is as well,
and Theorem 2.4.4 produces an isomorphism between
$HC^*_t(\ell^1(G))$ and the $\ell^1$-completion of the vector
space $\underset{x\in <G>}\to\oplus HC^*(\Bbb C)$.
\newpage

\vskip.2in
\centerline{\bf\underbar{
Appendix - The Baum-Connes assembly map and the generalized Bass Conjecture}}
\vskip.2in
\centerline{C. Ogle (OSU)}
\vskip.5in

In the early 1980's, P. Baum and A. Connes defined an assembly map
$$
\aA_*^{G,a}: KK^G_*(C(\underline{E}G),\Bbb C)\to K^t_*(C^*_r(G))
\tag{A.0.1}
$$
where $G$ denotes a locally compact group, $\underline{E}G$ the classifying space for proper $G$-actions, $C(\underline{E}G)$ the $G$-algebra of complex-valued functions on $\underline{E}G$ vanishing at infinity, and $KK^G_*(C(\underline{E}G),\Bbb C)$ the $G$-equivariant $KK$-groups of        
$(\underline{E}G)$ with coefficients in $\Bbb C$, while $K_*^t(C^*_r(G))$ represents the topological $K$-groups of the reduced $C^*$-algebra of $G$\footnote"$\dagger$"{For non-cocompact $\underline{E}G$, the equivariant $KK$-groups $KK^G_*(C(\underline{E}G),\Bbb C)$ are defined as a colimit over the directed set $Co_G(\underline{E}G)$ of $G$-cocompact subsets of $\underline{E}G$, with the relation given by inclusion:
$$
KK^G_*(C(\underline{E}G),\Bbb C):= \underset{Co_G(\underline{E}G)}\to{colim}KK^G_*(C(X),\Bbb C)
$$}. The original details of this map appeared (a few years later) in [BC1] and [BC2], with further elaborations in [BCH]. As shown in [BC3], when $G$ is discrete the left-hand side admits a Chern character which may be represented as
$$
ch_*^{BC}(G): KK^G_*(C(\underline{E}G),\Bbb C)\to \underset{x\in\ fin(<G>)}\to{\oplus} H_*(BG_x;\Bbb C)\otimes HPer_*(\Bbb C)
$$
where $fin(<G>)$ is the set of conjugacy classes of $G$ corresponding to elements of finite order, $G_x$ the centralizer of $g$ in $G$ where $x = <g>$, and $HPer_*(\Bbb C)$ the periodic cyclic homology of $\Bbb C$. Note that $H_*(BH;\Bbb C)\otimes HPer_*(\Bbb C)$ are simply the $2$-periodized complex homology groups of $BH$, and (via the classical Atiyah-Hirzebruch Chern character) can be alternatively viewed as the complexified $K$-homology groups of $BH$. Upon complexification, the map $ch_*^{BC}(G)$ is an isomorphism. The original construction of Baum and Connes $\aA_*^{G,a}$ was analytical. Motivated by the need to construct a homotopical analogue to their map, we constructed an assembly map in [Og1] which we will denote here as
$$
\aA_*^{G,h}\otimes\Bbb C: H_*(\underset{x\in\ fin(<G>)}\to\coprod BG_x;\tbrm{K}(\Bbb C))\otimes\Bbb C\to K_*^t(C^*_r(G))\otimes\Bbb C
$$
where $\tbrm{K}(\Bbb C)$ denotes the 2-periodic topological $K$-theory spectrum of $\Bbb C$. The construction of this map amounted to an extension of the classical assembly map constructed in [L] which was designed to take into account the contribution coming from the conjugacy classes of finite order. 

As before, let $\Bbb C[G]$ denote the complex group algebra topologized with the fine topology. To this algebra we associate the space 
$$
K^f(\Bbb C[G]) := \left|[q]\mapsto K_0(S_q(\Bbb C[G]))\times BGL^+(S_q(\Bbb C[G]))\right|
$$
where $\{S_q(X)\}$ denotes the total singular complex of the topological space $X$, and $|[q]\mapsto Y_q|$ the geometric realization of the simplicial set $Y$. This space has a natural infinite-loop structure coming from the standard (non-connective) Quillen-Grayson-Wagoner delooping of the spaces on the right, and the associated spectrum is a module over the spectrum delooping of $K^f(\Bbb C)$. In particular, there is  a Bott element   $\beta\in\pi_2(K^f(\Bbb C)$ (corresponding to the usual Bott element in $K_2^t(\Bbb C)$, and multiplication with this element induces a map
$$
\beta_*: \pi_*(K^f(\Bbb C[G]))\to \pi_{*+2}(K^f(\Bbb C[G]))
$$
For the purpose of this appendix, we define the \underbar{topological} $K$-theory of $\Bbb C[G]$ as
$$
K^t_*(\Bbb C[G]) := \pi_*(K^f(\Bbb C[G]))[\beta^{-1}]
$$
The two essential features of $\aA_*^{G,h}\otimes\Bbb C$, shown in [Og1], were (i) it factors through $K^t_*(\Bbb C[G])\otimes\Bbb C$, and (ii) the composition of $\aA_*^{G,h}\otimes\Bbb C$ with the complexified Chern-Connes-Karoubi-Tillmann character $ch^{CK}_*: K_*(\Bbb C[G])\otimes\Bbb C\to HC_*(\Bbb C[G])$ was effectively computable (see below). What we did not do in [Og1] was show that $\aA_*^{G,a}\otimes\Bbb C$ and $\aA_*^{G,h}\otimes\Bbb C$ agree. Since this initial work, there have been numerous extensions and reformulations of the Baum-Connes assembly map, as well as of the original Baum-Connes conjecture, which states that the map in (0.1) is an isomorphism. These extensions typically are included under the umbrella term \lq\lq Isomorphism Conjecture\rq\rq, (formulated for both algebraic and topological $K$-theory; cf. [DL], [FJ], [LR]). Thanks to [HP], we now know that the different formulations of these assembly maps (e.g., homotopy-theoretic vs. analytical) agree.
\vskip.2in

Abbreviating $KK^G_*(C(\underline{E}G),\Bbb C)$ as $K^G_*(\underline{E}G)$ (read: the equivariant $K$-homology of the proper $G$-space $\underline{E}G$), our main result is 

\proclaim{\bf\underbar{Theorem A.1}} There is a commuting diagram
$$
\diagram
K^G_*(\underline{E}G)\rto^{\aA_*^{G,DL}}\dto^{ch^{?}_*} & K^t_*(\Bbb C[G])\dto^{ch^{CK}_*}\\
HC_*^{fin}(\Bbb C[G])\rto|<<\tip & HC_*(\Bbb C[G])
\enddiagram
$$
where $\aA_*^{G,DL}$ is the homotopically defined assembly map of [DL], ${}^{fin}HC_*(\Bbb C[G]) :=  \underset{x\in\ fin(<G>)}\to{\oplus} HC_*(\Bbb C[G])_x\cong  \underset{x\in\ fin(<G>)}\to{\oplus} H(BG_x;\Bbb C)\otimes HC_*(\Bbb C)$ is the \underbar{elliptic summand} of $HC_*(\Bbb C[G])$, the lower horizontal map is the obvious inclusion, and the Chern character $ch^?_*$ becomes an isomorphism upon complexification for $*\ge 0$.
\endproclaim

Let $\beta$ denote a bounding class, $(G,L)$ a discrete group equipped with a word-length, and $H_{\beta,L}(G)$ the rapid decay algebra associated with this data. We write $K^t_*(H_{\beta,L}(G))$
for the Bott-periodic topological $K$-theory of the topological algebra $H_{\beta,L}(G)$. The \underbar{Baum-Connes} assembly map for $H_{\beta,L}(G)$ is defined to be the composition
$$
\aA_*^{G,\beta}: K^G_*(\underline{E}G)\overset{\aA_*^{G,DL}}\to\longrightarrow K^t_*(\Bbb C[G])\to K^t_*(H_{\beta,L}(G))
\tag{BC}
$$
where the second map is induced by the natural inclusion $\Bbb C[G]\hookrightarrow H_{\beta,L}(G)$.
In the introduction, we conjectured that the image of $ch_*:K^t(H_{\beta,L}(G))\to HC_*^t(H_{\beta,L}(G))$ lies in the elliptic summand ${}^{fin}HC_*^t(H_{\beta,L}(G))$ (conjecture $\beta$-SrBC). As the inclusion $\Bbb C[G]\hookrightarrow H_{\beta,L}(G)$ sends ${}^{fin}HC_*(\Bbb C[G])$ to ${}^{fin}HC_*^t(H_{\beta,L}(G))$, naturality of the Chern character $ch_*^{CK}$ and Theorem 1 implies

\proclaim{\bf\underbar{Corollary A.2}} If $\aA_*^{G,\beta}$ is rationally surjective, then $\beta$-SrBC is true.
\endproclaim

Since going down and then across is rationally injective, we also have (compare [O1])

\proclaim{\bf\underbar{Corollary A.3}} The assembly map $\aA_*^{G,DL}\otimes\Bbb Q$ is injective for all discrete groups $G$.
\endproclaim

We do not claim any great originality in this paper. In fact, Theorem A.1, although not officially appearing in print before this time,  has been a \lq\lq folk-theorem\rq\rq\  known to experts for many years. The connection between the Baum-Connes Conjecture (more precisely a then-hypothetical Baum-Connes-type Conjecture for $\Bbb C[G]$) and the stronger Bass Conjecture for $\Bbb C[G]$ discussed in the introduction was noted by the author in [O4]. 

There is some overlap of this paper with the results presented in [Ji]. A special case of Theorem 1 (for $* = 0$ and $\Bbb C[G]$ replaced by the $\ell^1$-algebra $\ell^1(G)$) appeared as the main result of [BCM].\vskip.3in

\head
Proof of Theorem A.1
\endhead

We use the notation $F_*^{fin}(\Bbb C[G])$ to denote the \underbar{elliptic summand} $\underset{x\in\ fin(<G>)}\to{\oplus}F_*(\Bbb C[G])_x$ of $F_*(\Bbb C[G])$ where $F_*(_-) = HH_*(_-), HN_*(_-), HC_*(_-)$ or $HPer(_-)$.
To maximize consistency with [LR], we write {\bf{S}} for the (unreduced) suspension spectrum of the zero-sphere $S^0$, {\bf{HN}}($R$) resp. {\bf{HH}}($R$) the Eilenberg-MacLane spectrum whose homotopy groups are the \underbar{negative cyclic} resp. \underbar{Hochschild} homology groups of the discrete ring $R$, and $\text{\bf{K}}^a(R)$ the non-connective algebraic $K$-theory spectrum of $R$, with $K_*^a(R)$ representing its homotopy groups. By [LR, diag. 1.6] there is a commuting diagram
$$
\diagram
H_*^G(\underline{E}G;{\text{\bf{S}}})\dto\rrto & & K^a_*(\Bbb Z[G])\dto^{NTr_*}\\
H_*^G(\underline{E}G;{\text{\bf{HN}}}(\Bbb Z))\dto\rto^{\cong} 
& HN_*^{fin}(\Bbb Z[G])\dto\rto|<<\tip & HN_*(\Bbb Z[G])\dto^{h_*} \\
H_*^G(\underline{E}G;{\text{\bf{HH}}}(\Bbb Z))\rto^{\cong} 
& HH_*^{fin}(\Bbb Z[G])\rto|<<\tip& HH_*(\Bbb Z[G])
\enddiagram
\tag{A.1.1}
$$
where the top horizontal map is the composition 
$$
H_*^G(\underline{E}G;{\text{\bf{S}}})\to H_*^G(\underline{E}G;{\text{\bf{K}}^a}(\Bbb Z))\overset{\aA^{G,DL}}\to\longrightarrow K_*(\Bbb Z[G])
$$
referred to as the the \underbar{restricted} assembly map for the algegraic $K$-groups of $\Bbb Z[G]$. The other two horizontal maps are the assembly maps for negative cyclic and Hochschild homology respectively. The upper left-hand map is induced by the map from the sphere spectrum to the Eilenberg-MacLane spectrum \text{\bf{HN}}, which may be expressed as the composition of spectra $\text{\bf{S}}\to \text{\bf{K}}^a(\Bbb Z)\to \text{\bf{HN}}$. By [LR], the composition on the left is a rational equivalence.

Let $\Bbb C^{\delta}$ denote the complex numbers $\Bbb C$ equipped with the discrete topology. Tensoring with $\Bbb C$ and combined with the inclusion of group algebras $\Bbb Z[G]\hookrightarrow \Bbb C^{\delta}[G]$, (A.1.1) yields the commuting diagram
$$
\diagram
H_*^G(\underline{E}G;\Bbb Q)\otimes\Bbb C\dto^{\cong}\rto & K^a_*(\Bbb C^{\delta}[G])\otimes\Bbb C\dto^{NTr_*}\\
HN_*^{fin}(\Bbb C[G])\rto|<<\tip & HN_*(\Bbb C[G])
\enddiagram
\tag{A.1.2}
$$
Next, we consider the transformation from algebraic to topological $K$-theory, induced by the map of group algebras $\Bbb C^{\delta}[G]\to \Bbb C[G]$ which is the identity on elements. By the results of [CK], [W] and [T], there is a commuting diagram
$$
\diagram
K_*^a(\Bbb C^{\delta}[G])\otimes\Bbb C\dto^{NTr_*}\rto & K^t_*(\Bbb C[G])\otimes\Bbb C\dto^{ch_*(\Bbb C[G])}\\
HN_*(\Bbb C[G])\rto & HPer_*(\Bbb C[G])
\enddiagram
\tag{A.1.3}
$$
where $ch_*(\Bbb C[G])$ is the Connes-Karoubi Chern character for the fine topological algebra $\Bbb C[G]$, and the bottom map is the transformation from negative cyclic to periodic cyclic homology. 

We can now consider our main diagram
$$

\diagram
H_*^G(\underline{E}G;\Bbb C)\otimes K_*(\Bbb C)\rto\dto & K_*^{a}(\Bbb C^{\delta}[G])\otimes\Bbb C\otimes K_*(\Bbb C)\rto\dto & K_*^{t}(\Bbb C[G])\otimes\Bbb C\otimes K_*(\Bbb C)\rto\dto^{ch_*(\Bbb C[G])\otimes ch_*(\Bbb C[\{id\}])} & K^t_*(\Bbb C[G])\otimes\Bbb C\dto^{ch_*(\Bbb C[G])}\\
HN_*^{fin}(\Bbb C[G])\otimes K_*(\Bbb C)\dto\rto & HN_*(\Bbb C[G])\otimes K_*(\Bbb C)\rto & HPer_*(\Bbb C[G])\otimes HPer_*(\Bbb C)\rto & HPer_*(\Bbb C[G])\\
HPer_*^{fin}(\Bbb C[G])\otimes K_*(\Bbb C)\rrto^{\cong}& & HPer_*^{fin}(\Bbb C[G])\otimes HPer_*(\Bbb C)\rto\uto|<<\tip & HPer_*^{fin}(\Bbb C[G])\uto|<<\tip
\enddiagram
\tag{A.1.4}
$$

The top left square commutes by (A.1.2), and the middle top square commutes by (A.1.3). The upper right square commutes by virtue of the fact that the Connes-Karoubi-Chern character is a homomorphism of graded modules, which maps the $K_*^t(\Bbb C)$-module $K_*^t(\Bbb C[G])$ to the $HPer_*(\Bbb C)$-module $HPer_*(\Bbb C[G])$, with the map of base rings induced by isomorphism $ch_*(\Bbb C[\{id\}]): K_*^t(\Bbb C)\otimes\Bbb C\overset\cong\to\longrightarrow HPer_*(\Bbb C)$. The lower left square commutes trivially, while the lower right commutes by the naturality of the inclusion $HPer_*^{fin}(\Bbb C[G])\hookrightarrow HPer_*(\Bbb C[G])$ with respect to the module structure over $HPer_*(\Bbb C)$.
Summarizing, we get a commuting diagram
$$
\diagram
H_*^G(\underline{E}G;\Bbb C)\otimes K_*(\Bbb C)\rto\dto^{\cong} &  K^t_*(\Bbb C[G])\otimes\Bbb C\dto^{ch_*(\Bbb C[G])}\\
HPer_*^{fin}(\Bbb C[G])\dto\rto|<<\tip & HPer_*(\Bbb C[G])\dto\\
HC_*^{fin}(\Bbb C[G])\rto|<<\tip & HC_*(\Bbb C[G])
\enddiagram
\tag1.5
$$
where the top map is the complexified Davis-L\"uck assembly map for topological group algebra $\Bbb C[G]$, and the bottom square is induced by the transformation $HPer_*(_-)\to HC_*(_-)$ (which respects the summand decomposition indexed on conjugacy classes). Restricted the elliptic summand yields the map$HPer_*^{fin}(\Bbb C[G])\to HC_*^{fin}(\Bbb C[G])$ which is an isomorphism for $*\ge 0$, implying the result stated in Theorem A.1.

\enddocument